\documentclass[11pt,regno]{amsart}

\usepackage{graphicx}
\usepackage[rgb]{xcolor}
\usepackage{amsmath,amscd,amsthm,amssymb,tikz} 
\usetikzlibrary{calc,positioning}
\usepackage[latin1]{inputenc}

\usepackage{overpic}

\usepackage{euscript,graphicx,color}
\usepackage{amsmath}
\usepackage{amssymb}
\usepackage{amsfonts}
\usepackage{mathrsfs}

\newtheorem{theo}{Theorem}[section]
\newtheorem{theo-app}{Theorem}[section]

\newtheorem{theorem}[theo]{Theorem} 
\newtheorem{corollary}[theo]{Corollary}
\newtheorem{proposition}[theo]{Proposition}

\newtheorem{lemma}[theo]{Lemma}
\theoremstyle{definition}
\newtheorem{example}[theo]{Example}
\newtheorem{remark}[theo]{Remark}
\newtheorem*{remark*}{Remark}
\newtheorem{definition}[theo]{Definition}
\newtheorem*{definition*}{Definition}

\newtheorem*{notation*}{Notation}

\theoremstyle{Theorem A}
\theoremstyle{Theorem B}
\theoremstyle{Theorem C}
\theoremstyle{Theorem D}
\theoremstyle{Theorem E}

\theoremstyle{Conjecture}

\newtheorem{notation}[theo]{Notation}

\theoremstyle{citing}

\numberwithin{equation}{section}

\renewcommand{\hat}{\widehat}

\DeclareMathOperator{\cl}{cl}

\newcommand{\R}{\mathbb{R}}

\newcommand{\C}{\mathbb{C}}
\newcommand{\Z}{\mathbb{Z}}

\DeclareMathOperator{\ind}{ind}

\DeclareMathOperator{\Per}{Per}
\DeclareMathOperator{\per}{per}

\DeclareMathOperator{\diam}{diam}

\DeclareMathOperator{\dist}{Dist}

\DeclareMathOperator{\Crit}{{\rm Crit}}

\DeclareMathOperator{\tree}{{\rm tree}}

\DeclareMathOperator{\reg}{{\rm reg}}

\DeclareMathOperator{\sym}{{\rm sym}}

\def\cC{\EuScript{C}}

\def\cM{\EuScript{M}}

\def\sR{\mathcal{R}}

 \def\sF{\mathcal{F}}

\newcommand{\cP}{\mathcal{P}}

            \def\d{\delta}

\def\ov{\overline}
\def\un{\underline}

 \def\dist{{\rm {dist}}}
          \def\s{\sigma}

\author[Feliks Przytycki ]{Feliks Przytycki} \address{Institute of Mathematics, Polish Academy of Sciences, ul. \'{S}niade\-ckich 8, 00-656 Warszawa, Poland}
\email{feliksp@impan.pl}

\begin{document}

\date{\today}

\title[Periodic orbits for polynomials]{There are not many periodic orbits in bunches for iteration of complex quadratic polynomials of one variable} 



\begin{abstract}

It is proved that for every complex quadratic polynomial $f$ with  Cremer's
fixed point $z_0$ (or periodic orbit)
for every $\delta>0$, there is at most one periodic orbit of minimal period $n$ for all $n$ large enough, entirely in the disc (ball) $B(z_0, \exp -\delta n)$
(and at most $p$ for a Cremer orbit of period $p$; in fact at most one, by a renormalization remark by C.L.Petersen). Next,  it is proved that the number of periodic orbits of period $n$ in a bunch $\cP_n$, that is for all $x,y\in \cP_n$, $|f^j(x)- f^j(y)|\le \exp -\delta n$ for all $j=0,...,n-1$, does not exceed $ \exp \delta n$.
We conclude that the geometric pressure defined with the use of periodic points coincides with the one defined with the use of preimages of an arbitrary typical point.  I. Binder, N. Makarov and S. Smirnov (Duke Math. J. 2003) proved this for all polynomials but assuming all periodic orbits were hyperbolic, and asked about general situations. We prove here a positive answer for all quadratic polynomials.

\end{abstract}


\maketitle

\tableofcontents

\section{Introduction}\label{Introduction}

In this paper we consider holomorphic maps $f:U\to \ov\C$, where $U\subset\ov\C$ is an open domain in the
Riemann sphere $\ov\C$, where $\C$ is the complex plane. Usually $f$ will be a rational function on $U=\ov\C$ or just a quadratic polynomial
$f(z)=f_c(z)=z^2+c$ for a complex number $c$.

\

The following notions of geometric pressure for a rational function $f:\ov\C\to\ov\C$ are of interest, see e.g. \cite{PRS2} and \cite{P-ICM18}. Let us start with just \emph{variational pressure} for continuous $f:X\to X$ for a compact set $X$ and a real-valued continuous function $\phi:X\to\R$ (called a \emph{potential}) function.

\medskip

 \begin{definition}[variational pressure]\label{var_pres}
\begin{equation}
P_{\rm var}(f,\phi)=
      \sup_{\mu\in{\cM}(f)}\left( h_\mu(f)+\int_X \phi \,d\mu\right),
      \end{equation}
where ${\cM}(f)$ is the set of all $f$-invariant Borel probability measures on $X$ and $h_\mu(f)$ is measure theoretical entropy. 
\end{definition}

For $f:\ov\C \to \ov\C$ a rational function with
and $\phi=-t\log |f'|$ we consider $f$ restricted to  $X=J(f)$ being Julia set for $f$, and any real $t$. Then we call the variational pressure above the \emph{variational geometric pressure} and $\phi$  the geometric potential. Though $\phi$ is not real continuous if $f$ has a critical point in $J(f)$, the definition still makes sense, since $\log|f'|$ is integrable for each $\mu$, \cite{P-Lyapunov}.


\begin{definition}[Tree pressure]\label{treep}
For every $z\in J(f)$ (or just $z\in \ov\C$) and $t\in \R$ define
\begin{equation}\label{treep-formula}
P_{\tree}(z,t)=\limsup_{n\to\infty}\frac1n\log\sum_{f^n(x)=z,\, x\in J(f)} |(f^n)'(x)|^{-t}.
\end{equation}
\end{definition}

\

These notions follow adequate notions in equilibrium statistical physics, e.g.
free energy for say generalizations of Ising model with a hamiltonian potential on the space of configurations
over $\Z^d$, see e.g. \cite{Ruelle} and \cite{Bowen}.

\

Another definition is also natural:

\begin{definition}[Periodic pressure]\label{periodicp}
For every $t\in \R$ define
\begin{equation}\label{periodicp-formula}
P_{\per}(t)=\limsup_{n\to\infty}\frac1n\log\sum_{f^n(z)=z,\, z\in J(f)} |(f^n)'(x)|^{-t}.
\end{equation}
\end{definition}

\

It was proved e.g. in \cite{PRS2}, see also \cite{P-conical}, that for every $z\in\ov\C$ except in a zero Hausdorff dimension subset of $\ov\C$
the tree pressure and variational geometric pressure coincide.

\medskip

It was proved in \cite[Theorem C]{PRS2}, that for $t\ge 0$ they coincide with the \emph{periodic pressure} for every rational function $f$ if the following Hypothesis $H$ holds,
see \cite{PRS2} and \cite{P-ICM18}.

Denote the standard spherical metric on the Riemann sphere $\ov\C$ by $\rho$. We assume
that the diameter of $\ov\C$ in $\rho$ is equal to 1.
For rational $f$ denote by $\Per_n$ the set of all periodic points of minimal period $n$ in the Julia set $J(f)$.

\

{\bf Hypothesis H}.\  For every rational function $f:\ov\C\to\ov\C$, for every $\d>0$ and all $n$ large enough, if for a set
$\cP\subset\Per_n$, for all
$x,y\in \cP$ and all $i:0\le i<n$ \
$\rho (f^i(x),f^i(y))< \exp -\d n$, then
$\# \cP\le \exp\d n$.

\

We prove here, in Section \ref{Hypothesis H},  that Hypothesis H holds for all quadratic polynomials.
Therefore we conclude with
\begin{theorem}\label{equality-of-pressures}
 For every quadratic polynomial
$f$ and real all $t\ge 0$
geometric pressures for $f$ acting in $\C$, variational and tree ones, are equal to the periodic pressure.
\end{theorem}

\

A condition stronger than Hypothesis H has been formulated in \cite{BMS}. Namely

{\bf Hypothesis BMS}. For every $\epsilon>0$ there exists $r=r(f,\epsilon)$ such that if $n\ge n(f,\epsilon)$ and $x\in\Per_n$, then
$$
\# \{y\in \Per_n : \rho (f^i(x),f^i(y))< r \forall 0\le i<n\} \le \exp \epsilon n.
$$
 The authors of \cite{BMS} proved it for all polynomials, but under the assumption REP saying that all periodic orbits of $f$ in $J(f)$ are hyperbolic repelling,
that is $|(f^n)'(x)|>1$ where $n$ is a period of $x$.
They used external rays, as we do here, and Yoccoz's puzzle structure (which we do not use).
They asked whether the assumption REP can be skipped.
See \cite[Lemma 7 and the comment preceding it]{BMS}.

\medskip

As we already mentioned, we succeed here, proving the equality of pressures in
Theorem \ref{equality-of-pressures} for all quadratic polynomials, without assuming REP, via proving  a weaker than Hypothesis BMS, but sufficient
for us, Hypothesis H.

\medskip

Hypothesis BMS, hence Theorem \ref{equality-of-pressures}, hold also for all generalized multimodal maps of interval, with all periodic orbits hyperbolic repelling, see \cite[Section 5]{PR}.

\

Let now  $f:U\to\ov\C$ be a holomorphic mapping as above. 
Let $z_0\in U$ be a fixed point for $f$. We call it indifferent if $|f'(z_0)|=1$. Then the  derivative $f'(z_0)$ is either a root of unity in which case we call $z_0$ parabolic, or not a root of unity when we call $z_0$ irrationally indifferent. In the latter case
if $f$ is not holomorphically linearizable in a neighbourhood of $z_0$ we call $z_0$ Cremer.
See e.g. \cite{Milnor-book} for an introduction to this theory. We use the same language for any periodic orbit of period $n$, replacing $f$ by $f^n$.

For any $x$ of period $n$ (in particular for $z_0$ as above) the derivative $(f^n)'(z_0)$ is sometimes called \emph{multiplier}.

\medskip

If $z_0$ is a  Cremer fixed point and $f'(z_0)=\exp \alpha 2\pi i$, where $\alpha$ is irrational then  its convergents $p_n/q_n$ in the continued fraction algorithm satisfy
\begin{equation}
\sum_{n\ge 1}\frac{\log q_{n+1}}{q_n} =\infty .
\end{equation}
(So convergence of the series, Bryuno condition, implies the existence of a so-called Siegel disc around $z_0$, that is a rotation in a disc by the angle $2\pi\alpha$ in respective holomorphic coordinates.)

If additionally
\begin{equation}
\sum_{n\ge 1}\frac{\log\log q_{n+1}}{q_n} <\infty
\end{equation}
then $f$ has a sequence of periodic orbits converging to $z_0$ of periods being a subsequence of $(q_n)$. see \cite{Perez-Marco1, Perez-Marco2}.

A question arises how many at most such periodic orbits converging to  $z_0$ in full, may happen.
In this paper we prove the following

\

\begin{theorem}\label{Cremer-fixed}
Set $f=f_c(z)=z^2+c$. Let $z=z_0\in \C$ be a fixed point for $f$. Then
 for every $\delta>0$ there exists $r_0>0$
 such that for every integer $n>1$
 and $r_n=r_0 \exp -n\delta$
 there is at most one orbit of minimal period $n$, entirely contained in $B(z_0,r_n)$.
\end{theorem}

\begin{theorem}\label{Cremer-periodic}
Set $f(z)=z^2+c$ as above. Let $z=z_0\in \C$ be a Cremer periodic point for $f$, of minimal period $p$.
Then
for every $\delta>0$ there exists $r_0>0$ such that for every $n>1$ and
 $r_n=r_0 \exp -n\delta$
 there are at most $p$ orbits of minimal period $np$, entirely contained in the union of the open discs $\bigcup_{i=0,...,p-1} B(f^i(z_0),r_n)$.
\end{theorem}
By estimates in Section \ref{distortion} the assumption "entirely" can be replaced by formally weaker
(in fact equivalent) assumption "at least one point of each orbit" in a disc (ball) or union
of the discs
in Theorems \ref{Cremer-fixed} or \ref{Cremer-periodic} respectively.

\medskip

Note that for every quadratic polynomial there is at most one Cremer periodic orbit. See \cite{Shishikura}.

\

A combination of these two theorems, a version not mentioning Cremer orbits, allows us to prove Hypothesis H, hence Theorem \ref{equality-of-pressures}.

\

An advantage when one deals with polynomials is a possibility to use external rays and angles (arguments), \cite{Milnor}, which makes the problem partially real, considering the angle  doubling map $F(\theta): = 2\theta \mod 2\pi \Z$, in place of complex. This was pioneered in \cite{BMS}. Here we apply Milnor's point of view of \emph{orbit portrait} for a Cremer periodic orbit $O(z_0)$ of period $p$, by rays starting at hyperbolic orbits of period $n$, close to $O(z_0)$, rather than at $O(z_0)$ itself.
For the reader's convenience we provide Milnor's approach in Appendix. We prove preservation of cyclic order, but we are not able to prove the unlinking property, which we fortunately can omit by a direct argument.

\medskip

We call sometimes the union of these rays, or various versions of truncated ones, \emph{quasi-spiders}, following a
terminology by J. Hubbard and D. Schleicher in \cite{Hub-Sch}, see Section \ref{Cremer_periodic}. Our
spiders are in fact multi-spiders because \emph{a priori} they correspond to finite families of periodic orbits,
rather than to individual periodic orbits.

\

Note that all our theorems proved here for quadratic polynomials hold for all unicritical polynomials, that is $z^d+c$ for $d\ge 2$, with the same proofs. In particular
the symmetric external rays, $\sR^{\sym}=-\sR$ in Subsection \ref{Step 3}, Figure \ref{Gamma}, can be replaced by the union of $d-1$ rays $e^{2\pi i b/d} \sR$ for $b=1,...,d-1$.

\medskip

A challenge is to verify our theorems for all polynomials of degree at least 2, where our approach might work.
For rational functions Theorem \ref{equality-of-pressures} holds if Hypothesis H holds. The latter holds
easily under some weak hyperbolicity assumptions, e.g. uniform exponential decay to 0 of diameters of pullbacks under $f^{-n}$ of small discs, see e.g. \cite{PRS1}. For a rational map the existence of a geometric coding tree, see \cite{P-ICM18}, without self-intersections, could make proving Hypothesis H doable, similarly to polynomials. But for general rational maps new ideas seem to be needed.

Of course our theory holds for any $g$ being a quadratic-like part of any rational $f$, see \cite{DH},  for its Julia set $J(g)$ because there is always a quasi-conformal (hybrid) conjugacy of $g$ to a quadratic polynomial, H\"older continuous, so preserving `small exponential' distances.
Compare Remark \ref{baby-M}.

\


{\bf Acknowledgement}. I am grateful to Genadi Levin for several helpful discussions. Supported by National Science Centre, Poland, Grant
OPUS 21 ``Holomorphic dynamics, fractals, thermodynamic formalism'',

\noindent 2021/41/B/ST1/00461.

\

\

\section{Distortion for iterates}\label{distortion}

\

We shall need the following auxiliary
\begin{lemma}\label{auxiliary}
Let $z_0$ be a fixed point for a holomorphic function $g:U\to \C$
for a domain $U\subset\C$, such that $|g'(z_0)|=1$.
Then, for every
$\alpha>0$ (in particular close to 0) there exists $n_0>0$ such that for
all $n\ge n_0$ and $i=0,...,n $ and every
$$
r\le r_n=r_0 \exp -n\delta
$$ as above,
it holds
\begin{equation}\label{inclusion}
g^i(B(z_0),r) \subset  B(z_0,r \exp  n \delta\alpha).
\end{equation}
\end{lemma}

\bigskip

\begin{proof}
Notice that except for a small neighbourhood of $\Crit(g)$, for $r_0$ small enough, for every 
$x\in B(z_0, r_0)$
\begin{equation}\label{Lip}
|g'(x)|\le 1 +L|x-z_0|,
\end{equation}
where
$$
L:=\sup_{x\in B(z_0, r_0)} \frac{ |g'(x)|-1}{|x-z_0|}.
$$
Let $n_\alpha$ be the smallest integer which satisfies

\begin{equation}\label{n-alpha}
1+ L r_0 \exp( -(n_\alpha+1) \delta\alpha)\le  \exp \delta\alpha.
\end{equation}
This roughly says that the exponential constant pace on the right hand side, dominates the parabolic pace  given by $|g'|$. Set $n'_\alpha: n_\alpha+1$

\medskip

Then, for each $j: 0\le j\le n-n'_\alpha$, 
we obtain

\begin{equation}\label{inclusion1}
g (B(z_0,   r\exp (j\delta\alpha)))
\subset B(z_0, r\exp ((j+1)\delta\alpha)).
\end{equation}

Indeed, denote $\eta:=r\exp (j\delta\alpha)))$. Then
\begin{equation}
\sup \{|g'(x)|: x\in B(z_0,\eta)\}
\le 1+2L\eta \le
\end{equation}
and using \eqref{n-alpha} we continue
\begin{equation}
\le 1+2L r_0\exp (-n\delta) \exp (n-n'_\alpha) \delta\alpha \le
1+2L r_0 \exp (-n'_\alpha )\delta\alpha)
\le \exp\delta\alpha.
\end{equation}
This indeed yields \eqref{inclusion1}

\smallskip

Now
compose \eqref{inclusion1},
starting from $j=0$ up to $j=i-1\le n-n'_\alpha$
for all $0\le i \le n-n_\alpha$. We obtain  

\begin{equation}\label{inclusion2}
g^i(B(z_0,r) \subset  B(z_0,r\exp ( 
i\delta\alpha)).
\end{equation}

Setting here $i=n-n_\alpha$, \underline{consider also
the last $n_\alpha$ iterates}. Denote the Lipschitz constant of $g$ by $M$.
Then

\noindent $g^{n}(B(z_0),r) 
\subset B(z_0, r\exp ((n-n_\alpha)\delta\alpha) M^{n_\alpha})\subset
B(z_0,r\exp (n\delta 2 \alpha))$
for $n$ large enough. The same for all $i\le n$.
We conclude that the inclusion \eqref{inclusion} holds for $2\alpha$
(which, if less than 1, can be considered as $\alpha$ in \eqref{inclusion}.

\end{proof}

The inclusion \eqref{inclusion} is self-reinforcing as follows:

\begin{corollary}\label{coro1}

For every $\delta>0$ and $\alpha:0<\alpha<1$
\begin{equation}\label{inclusion-inforced}
g^i(B(z_0,r_n)) \subset  B(z_0, r_n (1+ \lambda^{n}))
\end{equation}
for a constant $\lambda: 0<\lambda <1$,
depending on $\alpha$ and $\delta$,
for all $n$ large enough and $i=1,...,n$.
Moreover, the complex distortion of each $g^i$ on $B=B(z_0, r_n)$ is at most of order  $\lambda^n$, i.e.

\begin{equation}\label{complex distortion}
\sup_{x,y\in B} \Big\vert \frac{(g^i)'(x)}{(g^i)'(y)}-1\Big\vert\le \lambda^n.
\end{equation}
The same for every $r\le r_n$.

In particular the assertions hold for $\lambda$ arbitrarily close to $\exp (-\delta)$ for $\alpha$
close to 0.

Te assertions hold for all $n$ if $r_0$ is small enough.
\end{corollary}

\begin{proof}

For all $x\in B(z_0, r_n)$ for all $i\le n-n_{\alpha}$, due to \eqref{inclusion2} and \eqref{Lip}

\begin{align}\label{derivative}
|(g^i)'(x)| = & \exp\log \prod_{j=0}^{i-1} |g'(g^j(x))| \le
1+2\sum_{j=0}^{i-1}  (1+ L |g^j(x)-z_0)
\end{align}
$$ \le
 1+2\sum_{j=0}^{i-1} L r_n \exp ( j\delta\alpha)
 \le
 1+  2Lr_n \Bigl(\sum_0^{i-1} \exp j\delta\alpha\Bigr)
\le
1+ 2Lr_n \frac{(\exp i\delta\alpha) - 1}{(\exp \delta\alpha) - 1}.
$$

Hence, for  $C_{\delta,\alpha,L}:= 2L/(\exp \delta\alpha) - 1)$,
\begin{equation}\label{conclusion-inclusion}
g^i(B(z_0),r_n) \subset B(z_0,  r_n \sup_{x\in B(z_0,r_n)} |(g^i)'(x)|)
\end{equation}
$$
\subset B(z_0,  r_n (1+ C_{\delta,\alpha,L} r_0\exp (-n\delta)\exp (i\delta\alpha)).
$$
This includes all $i\le n$ if we increase adequately $C_{\delta,\alpha,L}$, as at the end of Proof of Lemma \ref{auxiliary}.
Thus, we prove Corollary \ref{coro1} with  any $1>\lambda > \exp (-\delta+\delta\alpha)$.
The inequality \eqref{complex distortion} follows from the calculation as for absolute values.
One can apply also Koebe Distortion Lemma, compare Remark below, to control the variation of the arguments of $(g^i)'$ by the variation of the absolute values.

Finally taking $r_0$ sufficiently small one gets  the assertions for first, finite number, of $n$'s immediately, having $g'-1$ arbitrarily close to 0.
\end{proof}

\begin{remark} 1. This Corollary, in particular \eqref{inclusion-inforced} does not need the holomorphy assumption. $C^2$-smoothness is enough. $|g'|-1$ is replaced by $||Dg(x)-Dg(z_0)||$ where $Dg(z_0)$ is the rotation operator.

2. In the holomorphic case a simple proof via Koebe Distortion Lemma is possible. Indeed,
consider the disc $B':=B(z_0, r_0\exp(-n\delta\alpha))$. Notice that proving \eqref{inclusion2} we could consider $B'$ in place of $B=B(z_0,r)$, due to iterated \eqref{inclusion1}. For $r_0$ small enough


each $g^i$ for each $i=1,...,n$, in particular $g^n$, is injective on $B'$ Consider the couple of discs
$$
B=B(z_0, r_n) \subset B'.
$$
Due to the injectivity of $g^n$ on the larger, right hand side disc, we can conclude, using Koebe Distortion Lemma, see e.g. \cite[Theorem 1.1]{Carleson-Gamelin} or \cite[Lemma 6.2.4]{PU}, that the distortion of $g^n$ on $B$, that is the fraction in \eqref{complex distortion},
 thats $(g^n)'(x)/(g^n)'(y)$ has absolute value bounded  by
$ 1+(8+\epsilon)x$, where
$x=\bigl(\exp (-n\delta )/\exp (-n\delta\alpha)\bigr) = \exp (-n(\delta-\alpha))$, for all $\epsilon>0$ and $n$ large enough.
Hence we are done with proving \eqref{complex distortion} by taking any
$\lambda: 1>\lambda> \exp (-(\delta-\alpha))$, in particular for any $\lambda>\exp (-\delta)$ if  $\alpha$ is taken appropriately.

\end{remark}

\section{Cremer fixed points. Proof of Theorem \ref{Cremer-fixed}}\label{Cremer}


\begin{proof}[Proof of Theorem \ref{Cremer-fixed}]

\

\subsection{ Preliminaries}\label{Step 1}

If $z_0$ is repelling or parabolic, the theorem follows from the topological behaviour
at its neighbourhood, see \cite[Camacho's theorem]{Milnor-book}. There are no periodic orbits there except $z_0$. This holds for
any holomorphic $f$. So we need consider only Cremer points.

In presence of a Cremer fixed point, Julia set of $f_c$ is
connected because it contains the only critical point 0 n it. Hence $A_c(\infty)$, the basin of attraction to $\infty$ by the action of $f_c$, is simply connected.
Let $\Phi_c: \{z: |z|>1\} \to A_c(\infty)$ be Riemann (holomorphic injective) map
from the outer disc in the Riemann sphere onto  $A_c(\infty)$, where $f(\infty)=\infty$ and $f'(\infty)=1$ (in the chart $z\mapsto 1/z$). Note that the function $\log |\Phi_c^{-1}|$ is Green's function $G=G_c$ on $A_c(\infty)$.

For every $\theta\in [0,2\pi\mod 2\pi \Z]$ we call the function
$\tau\mapsto\Phi (\tau \cdot (\exp  i\theta))$ for $1<\tau<\infty$ the external ray, with the external argument $\theta$.  We denote it or its range by $\sR_\theta$ and fixed  $\theta$  write $\sR(\tau)$. If there exists a limit
$\sR_\theta (\tau)$ as $\tau\searrow 1$ we say that the ray starts at this limit point.
We denote $\C$ compactified  by the circle at infinity $S(\infty)$ of external arguments, by $\widetilde\C$. We denote $f$ extended continuously to $\widetilde\C$ by $\widetilde f$.

\

Fix an arbitrary positive integer $n$.
Let $O_1,...,O_K$ be periodic orbits for $f_c$ of minimal period $n$, all in the disc $B(z_0,r)$ with $r\le r_n$. 
Each point $x$ in this family of orbits is a starting point of a nonempty finite number of external rays, $\sR(x,s)$.
The existence is known as Douady-Eremenko-Levin-Petersen theorem,
see \cite{P-accessibility} for generalizations.
Finiteness is a classical Douady-Hubbard easy observation.

 Namely $f$ permutes the family of the rays, hence external arguments, and the permutation
 restricted to external arguments for each $x\in O_k, k=1,...,K$ (leading to $f(x)\in O_k$) is the doubling map $F$ preserving cyclic order along
$ S(\infty)$, hence finite \cite[Definition preceding Theorem 1.1 and Lemma 2.3]{Milnor}. This preservation of the cyclic order of the rays starting (beginning) at each individual $x$
close to it (the order defined say via sectors between them) is clear, because $f$ is a homeomorphism there. It does not matter whether $x$ is periodic. Since the proof in \cite[Lemma 2.1]{Milnor} is only sketchy, we explain it below in [*] for completeness. Finally notice that this local order is the same as the order of the arguments (at $\infty$ ) by Jordan Theorem, see [*] and compare Subsection \ref{Step 2} below.
Note that $x$ need not be a fixed point for $f$ in these considerations; it is just in $J(f)$ and non-critical.

\medskip

[*] The order can be defined by taking any  $r: 0<r\le \infty$ (usually a small $r$) and considering for each $\sR=\sR(x,s)$ its truncation $\sR [r]$ joining the point $b\sR [r]$, being the last exit of $\sR$ from $\cl B(x,r)$, to $\infty$. Then the orders we consider are cyclic along the circles $S(x,r)=\partial B(x,r)$. Similarly along $S(\infty)$.

We prove that the order does not depend on $r$ by proving for each $r<\infty$ that it coincides with  the order in $S(\infty)$.
Indeed, consider an arbitrary $r<\infty$.
Consider three rays $\sR_j$ beginning at $x$ and their truncations $\sR_j[r]$ for $j=1,2,3$,
$b\sR_3 [r] \in S(x,r)$ belonging to the arc $\gamma_{1,2}[r]$ between $b\sR_1[r]$
and $b\sR_2 [r]$
counterclockwise.
Then by Jordan theorem, $\sR_3 [r]$ must be in the Jordan domain $Q$ in $\widetilde\C$ bounded by $\sR_1[r], \sR_2[r]$,
$\gamma_{1,2}[r]$ and $\gamma_{1,2}(\infty)\subset S(\infty)$, the latter between the arguments of $\sR_1$ and $\sR_2$ counterclockwise. Since $\sR_3[r]$ stays in $Q$ and tends to $\infty$, it must tend to a point (argument)
in $\gamma_{1,2}(\infty)$. Thus the orders in $S(r)$ and $S(\infty)$ coincide.

\

[**] Notice that we can change this definition by replacing the last exit by the first exit, making the definition  purely local.
Indeed, consider $0<r_1<r_2$ and denote the first exits of $\sR_j[r_i], i=1,2$ from
$ B(x,r_i)$ (Namely the first point in $\C \setminus B(x,r_i)$) by $e\sR_j[r_i]$ and the sub-ray of $\sR_j[r_i]$ joining $x$ to $e\sR_j[r_i]$ by
$\sR_j[x,r_i]$. Denote the arc in $S(x,r_i)$ joining $e\sR_1[r_i]$ and $e\sR_2[r_i]$ counterclockwise, by $\gamma_{1,2}[r_i]$. Denote  Jordan domains bounded by $\sR_1[x,r_i], \sR_2[x,r_i]$ and $\gamma_{1,2}[r_i]$  by $U[r_i]$ for $i=1,2$.

Assume that $e\sR_3[r_1]\in \gamma_{1,2}[r_1]$. Then $\sR_3[r_1]\subset U[x,r_1]$ except its ends.
 So its part close to $x$ is in $U[x,r_2]$, because $\sR_j[2r] \setminus\sR_j[r], j=1,2$ and $\gamma_{1,2}[r_2]$ run within a positive distance from $x$.
So $\sR_3$ stays in $U[x,r_2]$ until it leaves it, necessarily through $\gamma_{1,2}[r_2]$  since the rays cannot intersect one another.
Thus the orders in the circles
$S(x,r_1)$ and $S(x,r_2)$ coincide.

Notice finally that the orders via the "first" and the "last" exit coincide, because the "last" from $x$ is the "first" from $\infty$.

\medskip

We can replace above any circle by a loop around $x$ or $f(x)$. So $f$ preserves the local cyclic order,
defined with the use of the loops $S(x,r)$ and $f(S(x,r))$. See also Subsection \ref{Step 2}.

\medskip

But we shall prove more, that the preservation 
by $f$ of this cyclic order holds for the set
of \underline{all external arguments} of the family $\bigcup_x\sR(x,s)$
(as if all $\sR(x,s)$ started at $z_0$). See Corollary \ref{coro-order}.

\noindent \underline{This is the main point of the paper.}.

\

\subsection{Sub-rays and sectors}\label{Step 2}

 \

For an arbitrary $r: r_n \le r\le 2r_n$, consider the part $\sR'(x,s)=\sR'(x,s)[r]$ of  $\sR(x,s)$
starting at \underline{the last intersection with the circle}

\noindent \un{$S(r):=\partial B(z_0,r)$} and going to $\infty$. Enumerate these rays along the cyclic order in $S(r)$, writing 
$\sR'_j[r]$ or $\sR_j[r]$, or just $\sR'_j$ if $r$ is fixed. The non-truncated rays $\sR(x,s)$ will be denoted by $\sR_j$.

Here $m=Nn$, where $n$ is the minimal period, as before, and $N=Kv$, where $K$ is the number of periodic orbits as above, and $v$ is the valence, namely the number of $\sR(x,s)$ starting at each $x$.
\emph{A posteriori} it will be clear that $v$ does not depend on $x$. In fact $v=1$ or $2$ for $n>1$, by  Proposition \ref{valence}. So $Kn\le m\le 2Kn$. We will get even  $K=1$ by Theorem \ref{Cremer-fixed}.

We shall consider the domains,
\emph{\bf sectors}, $Q_j[r]\subset \C, j=0,...,m-1$
bounded respectively by $\partial Q_j[r]$ consisting of
$\sR'_j$, $\sR'_{j+1 (\mod m)}$ and the arc $\Gamma_j[r]$ in the circle joining their starting points, not containing others.
We can consider $\ov Q_j[r]$ and $\widetilde Q_j[r]$, the closures of $Q_j[r]$ in $\ov\C$ or in $\widetilde\C$
respectively.

\medskip

Note that the cyclic order above coincides with the cyclic order of the arguments of $\sR(x,s)$'s
at $\infty$.
The same holds if $S(r)$ is replaced by any small simple loop around $z_0$ and our family of $O_k$'s.

The proof is the same as in [*] where all the rays begin at one point, but we repeat it for completeness. Indeed, no $\sR'_i[r]$ starting  $S(r)\setminus \Gamma_j[r]$, namely outside the arc
in $S(r)$ joining beginnings of the rays $\sR'_j[r]$ and $\sR'_{j+1 \mod m}[r]$, can enter the sector $Q_j[r]$ to go to $\infty$ in it, by Jordan Theorem. It cannot intersect $\sR'_j[r]$ nor $\sR'_{j+1}[r]$ because the rays do not intersect each other. Nor it intersects $\Gamma_j[r]$ because it would intersect $B(z_0,r)$ before, but its beginning is the last point in $\cl B(z_0,r)$ by definition.
We conclude that the order does not depend on $r$ (nor on any loop above).

\medskip

Denote the permutation
induced  by $f$ on the set of the rays, in the other words on the indices $\{0,1,...,m-1\}$,
by $\sigma$.

\medskip

 For all $\sR'_j$ we can find sub-rays $\widehat\sR_j\subset\sR'_j$, ending at $\infty$ and beginning
$b\widehat\sR_j$ in the disc $B(z_0,2r)$, that is

 \begin{equation}\label{beginnings}
 r\le |b\widehat\sR_j-z_0|<2r
 \end{equation}

such that

 \begin{equation}\label{invariance}
 f(\widehat\sR_j)\subset \widehat\sR_{\sigma(j)} \subset \C\setminus B(z_0,r)
  \end{equation}
  for each $j$. 

Indeed, for each $\sR'_j$ we find $j_{\max}$ such that in its $f$-orbit
the beginning of $\sR'_{j_{\max}}$ has maximal value of Green's function $G_c$ (see the beginning of Subsection \ref{Step 1}.
Define for $j=0,...,m-1$
$$
\widehat\sR_{\sigma^j(j_{\max})}=f^{j-j_{\max}}(\sR'_{j_{\max}}).
$$
Since under the action by $f$ the value of
$G$ grows by $\log 2$, the inclusion \eqref{invariance} holds (it is equality for
$j=0,...,m-2$).
Moreover  $\widehat\sR_j\subset \sR'_j$ holds.
Finally beginnings of all $\widehat\sR_j$ are in $B(z_0,2r)$ by Corollary \ref{coro1}.

Though Lemma \ref{auxiliary} and Corollary \ref{coro1} hold for $0<i<n$ and here we need it for $i\le Nn$,
we can just write $r_n=r_0\exp -n\delta= r_0\exp -Nn (\delta/N)$ ( so that $\delta/N$ takes over the role of $\delta$),  and correct $r_0$ appropriately.

\medskip

It will be convenient to replace each $\widehat\sR_j$ by its sub-ray
(by \eqref{beginnings})  $\sR'_j[2r]$.

Summarizing
\begin{equation}\label{summarizing}
\sR'_j[r] \supset \widehat\sR_j[r] \supset\sR'_j[2r].
\end{equation}

With some more care (assuming $r$ sufficiently small) we can assure that the latter truncated rays
satisfy
\begin{equation}\label{32}
\sR'_j[\frac32 r] \supset\widehat\sR_j[\frac32 r] \supset \sR'_j[2r].
\end{equation}
To see it, we define here $\widehat\sR_j$ as above, but starting from $S(\frac32 r)$ in place of $S(r)$.
We conclude with

\begin{equation}\label{32c}
f(\sR'_j[2r]) \subset \sR'_{\sigma(j)}[\frac32 r].
\end{equation}

\

While considering truncated rays $\sR'$ where the distance from $z_0$ of their beginning is specified
we omit `prime' writing e.g. $\sR[2r]$ in place of $\sR'[2r]$ as above, see the beginning of subsection \ref{Step 2}.

 \

Note that the union $\bigcup_j \partial Q_j[2r]$ has no self-intersections (as for any other radius, not just $2r$).

Note that despite \eqref{summarizing}, it
might not be true that $Q_j[r] \supset Q_j[2r]$, see Figures 
\ref{backward} and \ref{Gamma}.

\


\begin{definition}[critical and regular sectors]\label{regular-critical}
Assume in our notation that $Q_0[2r]$ is the only sector, called the \emph{critical sector},  containing the critical point 0. Denote the union of all other sectors (together with the rays separating them), thus called the \emph{regular sector}, by $Q^{\reg}[2r]$.
Each individual sector in the union
will also be called a \emph{regular sector} and its boundary arc in $S(\infty)$ a \emph{regular boundary  arc}.

Denote $\Gamma^{\reg}[2r]:=\bigcup_{j=1,...,m-1} \Gamma_j [2r]$.
Similar objects will be considered for $r$ or any (small) $r'\ge r$ if we fix $r$ from now on.

Also for any $0\le j_1< j_2 \le m$ define
$Q_{j_1,j_2}[r']:=\bigcup_{j=j_1,...,j_2-1} Q_j[r']$
\end{definition}

\medskip

The following holds for $F$ being extension of $f$ to $S(\infty)$ 

\begin{proposition}\label{regular-infty-arc} The arc
$\partial_\infty$ in $S(\infty)$ bounding the regular sector $Q^{\reg}[2r]$, does not depend on $r$
for $r$ small enough. 
\end{proposition}

{\bf Remarks} 1.
This fact follows a posteriori from the Main Lemma  \ref{injectivity} in Subsection \ref{Step 6}, which yields that this arc is the complement in $S(\infty)$ of the only arc between the consecutive external arguments of $\sR_j$ and $\sR_{j+1 \mod m}$ longer counter-clock-wise than $\pi$.

\smallskip

2.  Denote by $k_1, k_2$ the indices $j$ at $\sR_j$ bounding $Q^{\reg}[2r]$. Since a priori they
may depend on $2r$, they would be denoted $k_1[2r], k_2[2r]$, and only a posteriori we learn that they do not. See Proposition
\ref{regular-infty-arc}. However, for simplification, fixed $2r$ we shall omit writing it at $k$.

\smallskip

3.  At Figure \ref{backward} the critical point $0\notin Q^{\reg}[2r]$ but it  belongs to the sector $Q_{k_1[2r],k_2[2r]}[r]$, that is with the same $k_1,k_2$ but $2r$ replaced by $r$ to mark
the bottom arc of the sector (called the \emph{floor} later on).
So the arcs bounding regular sectors for $2r$ and $r$ are different, contradicting
Proposition \ref{regular-infty-arc}. So the picture is impossible.

A direct argument is shown at Figure \ref{twin}. The picture shows that the growth of the arguments $\Delta_0 \sR_k[r,2r]$ and $\Delta_c f(\sR_k[r,2r])$ with respect to 0 or $c$ respectively
are both approximately $2\pi$ which is impossible because 0 is a critical point so the latter
should be around $4\pi$.

\smallskip

4. Here and in the sequel we write $\sR_k[r,2r]:=\sR_k[ r] \setminus \sR_k[2r]$  for $k=k_1,k_2$.
Remember that we consider here $k=k[2r]$.
In particular, for $r'\ge r$, writing $Q_0[r']$ and $Q^{\reg}[r']$, we mean $k=k[2r]$
Denote also the beginning and end points of
$\sR_k[r,2r]$ by $b\sR_k[r,2r]$ and $e\sR_k[r,2r]$ respectively.

 \begin{figure}

\includegraphics[height=6cm,width=12cm,trim=100pt 100pt 100pt 100pt,clip]{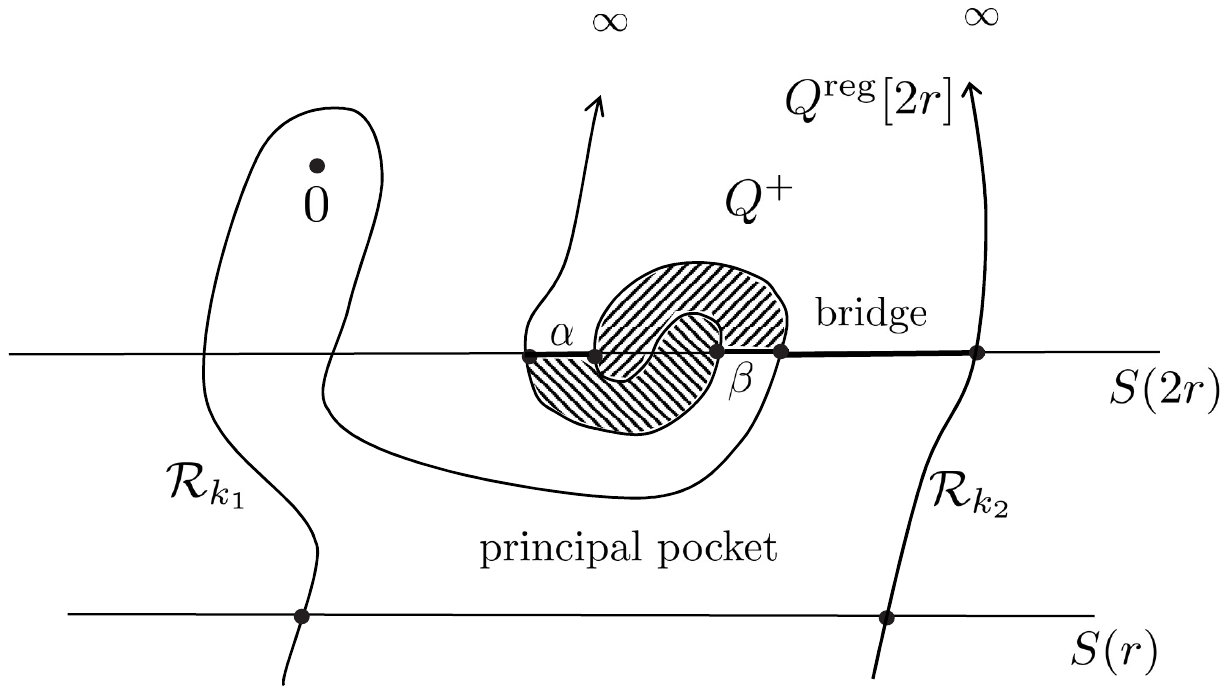}

\caption{Here $\alpha$ is an upper arc and $\beta$ is a lower arc. The adjacent pockets are hatched NE and NW respectively.} 
\label{backward}
\end{figure}

\begin{figure}[h]

\includegraphics[height=9cm, width=12cm, ,trim=50pt 120pt 0pt 100pt,clip]{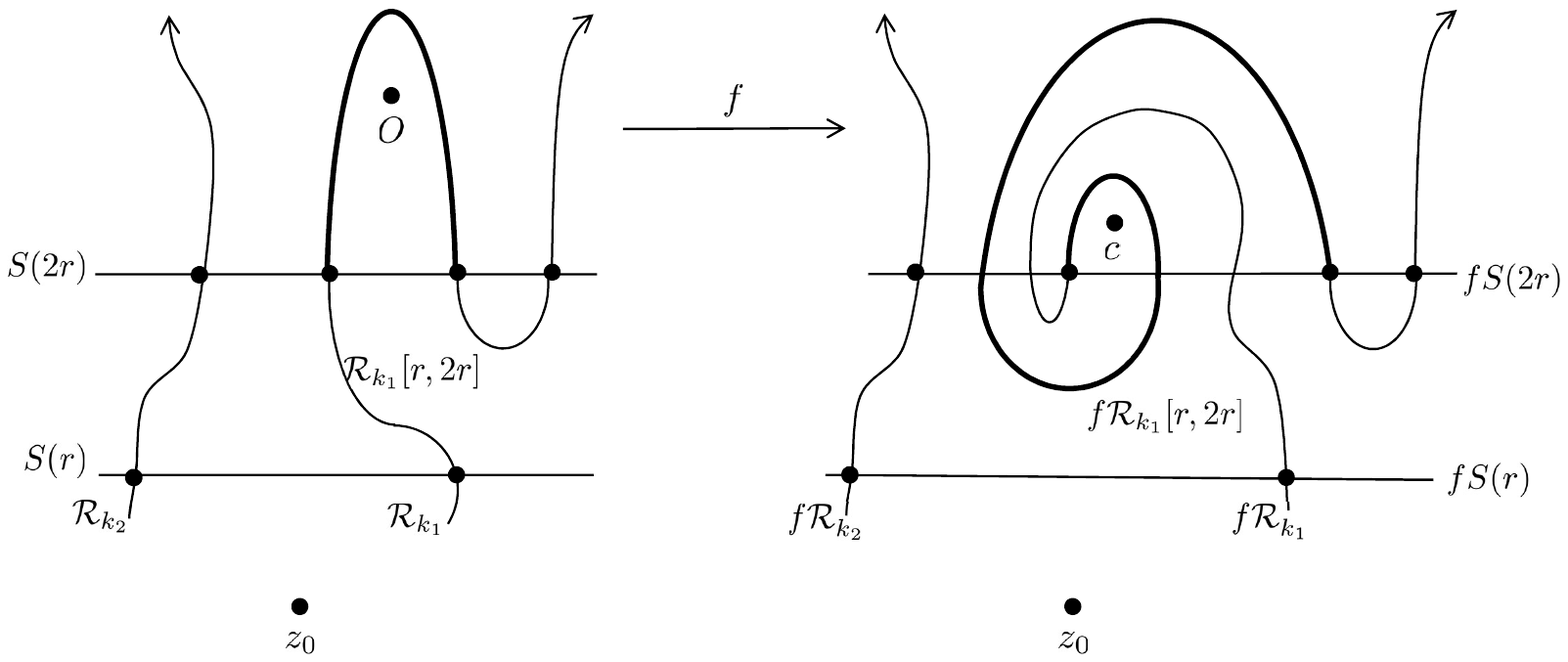}

\caption{An explanation why Figure \ref{backward}, with $0\in Q^{\reg}[r]$ and
$0\in Q^{\reg}[r]$, is impossible.}
\label{twin}
\end{figure}

\smallskip

5. Our further considerations will hold in fact for any $1\le k_1[2r]<k_2[2r]\le m$, so such that
the critical point 0 is not in the related sector.

\

\underline{Until the end of this Subsection we will be introducing a basic}

\noindent \underline{terminology}, to be used in this paper only partially, but substantial in future we believe.

\medskip


\

\

{\bf Terminology and Remarks}

\medskip

For any $K\in \sR_k[r,2r], k=k_1,k_2$ denote by $\sR_k[K,2r]$ the sub-curve of $\sR_k[r,2r] $
beginning at $K$ ending at $e\sR_k[r,2r]$.

\bigskip

1.  Consider three curves: $\Gamma^{\reg}[2r]$ and $\sR_k[r], k=k_1,k_2$. Their union in $Q^{\reg}[r]$
cuts it into several components. The one
 reaching $\infty$ will be called the \emph{unbounded upper domain}. Denote it  by $Q^+$.

 \medskip

The boundary in $\C$, $\partial_\C Q^+$, contains exactly one arc in $\Gamma^{\reg}[2r]\cap \partial Q^+$ joining  $\sR_{k_1}[r,2r]$ with $\sR_{k_2}[r,2r]$, intersecting these rays only at its ends. Call it sometimes
 an \emph{upper bridge} or just a \emph{bridge}.

Indeed, fix a convention that $e\sR_{k_1}[r,2r]$ is "to the left" of $e\sR_{k_2}[r,2r]$ in
$\Gamma^{\reg}[2r]$. Denote the most to the right point of $\sR_{k_1}[r,2r]\cap \Gamma^{\reg}[2r]$ by
$\widehat K^{k_1}$ and the most to the left point of $\sR_{k_2}[r,2r]\cap \Gamma^{\reg}[2r]$ by
$\widehat K^{k_2}$. If  $\widehat K^{k_2}$ is to the right of $\widehat K^{k_1}$, then we define
our bridge as $(\widehat K^{k_1}, \widehat K^{k_2})$. If $\widehat K^{k_2}$ is to the left of $\widehat K^{k_1}$, then we have two cases:

I. $\widehat K^{k_1}$ is "below" $\sR_{k_2}[\widehat K^{k_2},2r]$ to the right of $\widehat K^{k_2}$, that is
it is not in $\partial Q^+$. See Figure \ref{channel}. Then define
${K^\circ}^{k_1}$ as the most to the right point in $\sR_{k_1}[r,2r]\cap \Gamma^{\reg}[2r]$
to the left of $\widehat K^{k_2}$. Then define our bridge as $({K^\circ}^{k_1}, \widehat K^{k_2})$.

II. $\widehat K^{k_2}$ is "below" $\sR_{k_1}[\widehat K^{k_1},2r]$ to the left of $\widehat K^{k_1}$.
Then, swapping above $k_1$ with $k_2$ we define ${K^\circ}^{k_2}$ and the bridge is
defined as $(\widehat K^{k_1}, {K^\circ}^{k_2})$.

\medskip

The existence of the bridge and its uniqueness follow immediately from the construction.

\smallskip

\begin{figure}[h]

\includegraphics[height=7cm, width=11cm, trim=0pt 50pt 0pt 100pt,clip]{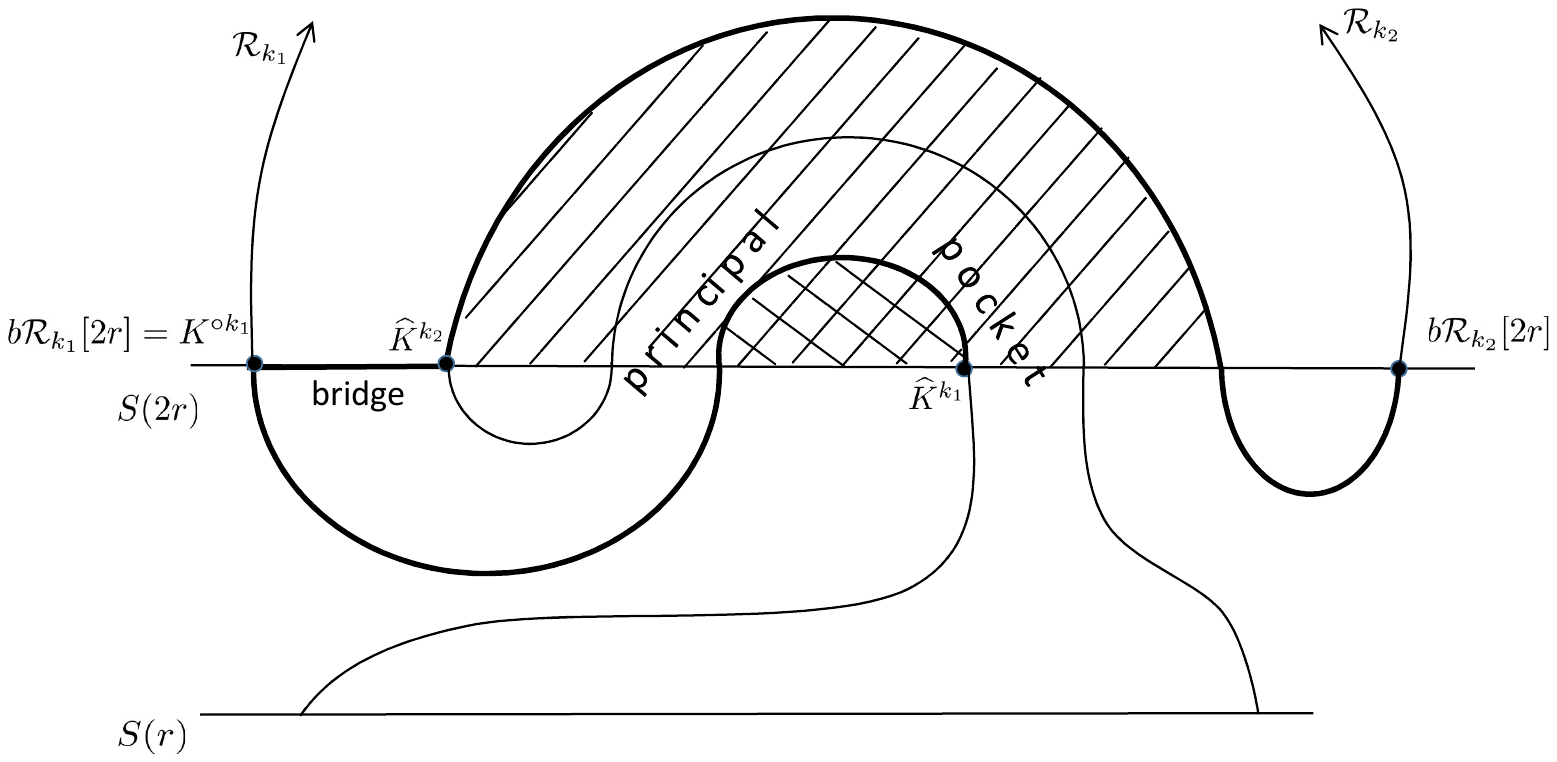}

\caption{Construction of the bridge and principal pocket.
A lower pocket for $k_2$ is hatched NE. Cross-hatched is a lower pocket for $k_1$.}
\label{channel}
\end{figure}

\smallskip

 2. The sub-arcs of $\Gamma^{\reg}[2r]$ in $\partial Q^+$ with both end points either in $\sR_{k_1}[r,2r]$ or both in $\sR_{k_2}[r,2r]$,
 are called the \emph{upper arcs}. Denote their end points by $P$ and $P'$ and the
 arcs by $(P,P')$ (warning: $K,K'$ and $L,L'$ to be defined later on, see \eqref{subarcs} are different).
 Call the Jordan domain bounded by $(P,P')$ and the
 sub-curve $\Upsilon_P\subset \sR_k[r,2r]$ joining the ends $P,P'$ of the arc, an \emph{upper  pocket}.

3.  Finally
consider for each $k=k_1,k_2$ the curve
$\sR_k[\widehat K^k,2r]$ in $\sR_k[r,2r]$ defined above joining $\widehat K^k$ to $e\sR_k[r,2r]$,
and consider the components to which it, together with $\Gamma=\Gamma^{\reg}[2r]$,
cut $\C\setminus Q^+$. Denote by $Q_k^-$ the only unbounded component.
Among the bounded components consider those whose borders share an arc with
$\partial Q_k^-$.
 Such an arc is a sub-arc of $\Gamma^{\reg}[2r]$ joining two points in $\sR_k[\widehat K^k,2r]$, with the same $k$. Call this arc a \emph{lower arc},
 and the component bounded by it and the curve $\Upsilon_P\subset \sR_k[\widehat K^k,2r]$ joining the ends $P,P'$ of the arc, a \emph{lower pocket}, by analogy to the upper arc and the upper pocket, but on the "bottom side" of
 $\sR_k[\widehat K^k,2r]$,
 unlike the "upper side" before. Notice that the pockets for $k_1$ need not be disjoint from the pockets for $k_2$, see Figure \ref{channel}.

\smallskip

4. 
The Jordan domain bounded by the bridge, the subcurves of $\sR_k[r,2r]$ joining
beginnings of $\sR_k[r]$ with the ends of the bridge
and $\Gamma^{\reg}[r]$
is called the \emph{principal pocket}.

\medskip

5.  The arc $\beta$ (\emph{lower arc}) on Figure \ref{backward} enables for its $f$-image, intersecting of
 $f(\sR_k[2r])$ with $\Gamma'$ and needs care.
 Moreover it is not excluded a non-empty  intersection of $f(\sR_{k_1}[2r])$ with a lower pocket with a sub-curve of $f(\sR_{k_2}[\widehat K^{k_2},2r])$ in its boundary  , or a similar intersection with $k_1$ and
 $k_2$ interchanged. See Figure \ref{channel}.

 On the other side $f(\alpha)$ and $f(\sR_k[2r])$
   cannot intersect one another
 by Lemma \ref{no-upper}. They could, if $f$ had at least 2 critical values in two neighbouring
critical sectors, to be elaborated in a future paper by the author.


\

\subsection{The Fold Lemma and homoclinic-like loops}\label{Step-Fold}

\

\

The following fact will be useful

\begin{lemma}[Fold Lemma]\label{fold}

Let $f=f_c$ be a quadratic polynomial. Let $\gamma_1, \gamma_2$ be two continuous curves in $\C$ such that
the beginnings and ends of their $f$-images coincide, that is
$f(b\gamma_1)=f(b\gamma_2)$ and $f(e\gamma_1)=f(e\gamma_2)$.
Suppose that $\delta:=f(\gamma_1) \cup f(\gamma_2)$ does not contain the critical value $c$. Suppose $\delta$ is a Jordan curve, in particular that it bounds a simply-connected domain $Q$. Then

1. $c\in Q$ if and only if $b\gamma_1=b\gamma_2$ and $e\gamma_1)\not=e\gamma_2$ (then call
$f(b(\gamma_i))$ a vertex and $f(e(\gamma_i))$ a homoclinic-like point), or
the same with changed the roles of $b$ and $e$.

2. $c\notin Q$ if and only if $b\gamma_1=b\gamma_2$ and $e\gamma_1=e\gamma_2$
or $b\gamma_1\not=b\gamma_2$ and $e\gamma_1\not=e\gamma_2$.

\end{lemma}

\begin{proof}
(classical, provided for completeness.)

Indeed, if $c\notin Q$, then fixed any point $a_0\in Q$ and $a_1\in f^{-1}(a_0)$ we can define $g_1(a_0)=a_1$ and extend holomorphically $g$ along any curve in $Q$. This defines
$g_1$ at $z$ the end of the curve, not depending  on the choice of the curve joining $a_0$
with $z$ by the monodromy theorem, even for $z$ in a neighbourood of $\cl Q$. $g_1$ is a branch of $f^{-1}$ and $g_1(Q)$ is a Jordan domain $Q_1$. The other choice $g_2(a_0)=a_2\in f^{-1}(a_0)$ yields another branch $g_2$
of $f^{-1}$ on  $Q$ and $Q_2=g_2(Q)$ is Jordan. This construction extends to a neighbourhood  $\cl Q$ so $\cl Q_1$ and $\cl Q_2$ are disjoint. This implies the property described in the item 2.

\smallskip
\
Since $f$ is of degree 2, $f^{-1}(\delta)$ consists of one or two disjoint Jordan  curves $\delta_1$ and $\delta_2$. If $c\in Q$ then in the latter case $f:\delta_i \to\delta$
has degree 1 for $i=1,2$, so the critical point 0 is in none of the Jordan domains $V_i$ bounded by $\delta_i$ and by maximum principle $f$ is a homeomorphism to $Q$,
so $c\notin Q$,
a contradiction. So $\delta':=f^{-1}(\delta)$ is one Jordan curve and degree of
$f:\delta'\to \delta$ is 2. This yields the situation described in the item 1.

\end{proof}

Note that this Lemma holds for every continuous orientation preserving branched covering of degree 2 of $\C$.

\smallskip

\begin{remark}\label{subfolds}
 The assumption that $f(\gamma_1)$ and $f(\gamma_2)$ do not intersect, except at their beginning and end points, can be omitted in the following sense: Consider any two different intersection points
$f(a_1)=f(a_2)$ and $f(b_1)=f(b_2)$ for $a_1,b_1\in \gamma_1$ and $a_i,b_i\in \gamma_i$
for $i=1,2$ the $f$-images closest to each other in one of $f(\gamma_i)$. Consider the curves $\gamma'_i$ joining $a_i$ with $b_i$ in $\gamma_i$, for $i=1,2$.
Then
the assertion of Lemma \ref{fold} holds for the Jordan curve $\delta':=f(\gamma'_1) \cup f(\gamma'_2)$ bounding a Jordan domain $Q'$.
\end{remark}

\smallskip

In applications one curve will be an arc in a circle, say $S(2r)$, another a sub-curve in an
external ray, say in $\sR_k[2r]$. We can call its $f$-image (with respect to the $f$-image of the arc, a \emph{fold}, yielding the name of the lemma. We could
call the pair $f(\gamma_1), f(\gamma_2)$ a  \emph{ homoclinic-like loop}, bounding a Jordan domain $Q$. If $Q\ni c$ we call it a \emph{critical loop}. If not, we call it
a \emph{non-critical loop}.

\

\subsection {On the rays bounding the regular sector.
Modifying the "floor"}\label{Step 3}

\

\

Consider  $f(\partial Q^{\reg}[2r])\subset\C$. If this curve does not have any self-intersections, the proof of the preservation of the cyclic order as at the end of Subsection \ref{Step 1} is easy, see Subsection \ref{Step 6}, especially Main Lemma \ref{injectivity}.  The non-existence
 of the self-intersections clearly holds if $z_0$ is repelling and all the rays under consideration begin at it,
and $r$ is adequately small.

In our situation
the only self-intersection points which  $f(\partial Q^{\reg}[2r])$ could have, are the points of
$f(\Gamma^{\reg}[2r])
\cap f(\sR_k [2r])$, for $k=k_1,k_2$.

Call the arc $\Gamma:=\Gamma^{\reg}[2r]$
the \emph{floor of the regular sector}
We shall remove the self-intersections by replacing the floor, hence its $f$-image
$$
\Gamma':=f(\Gamma)\subset f(\partial Q^{\reg}[2r])
$$
by a homotopic curve, $(\Gamma^*)':=f(\Gamma^*)$, the homotopy in $\C\setminus \{c\}$, in particular with the same growth argument with respect to $c$.

In fact we shall replace $\Gamma$ by a homotopic curve $\Gamma^*$ the ends intact, disjoint from $f^{-1}(f(\sR_k[2r]))$, so that the image of the homotopy omits $c$.

Thus, we shall prove that the growth
 of the argument along the curve $f (\partial Q^{\reg}[2r])$ in $\C$ with respect to the critical value $c$ is as if the self-intersection does not happen
 at all, unlike on Figure \ref{impossible}, on the left hand side.

We must be careful, because the rays seem to be "almost dense" in the Julia set $J(f)$ and it is easy to see that if the self-intersections happen for all small $r$ then $J(f)$ is not locally connected (because $f(\sR[2r])$ must then pass
arbitrarily close to both: $z_0$ and its pre-image $-z_0$, and close to 0.

Indeed if a polynomial has a Cremer fixed (or periodic) point, then its Julia set is not locally connected, see \cite[Theorem 18.5]{Milnor-book}.

\

The following lemma below simplifies the construction of $\Gamma^*$.

 \medskip

\begin{notation}[The branches of $f^{-1}$]\label{not:branches}
Denote  by $g$ the local branch of $f^{-1}$ mapping $z_0$ to itself.
Then $g(\Gamma')=\Gamma^{\reg}[2r])$. The second local branch, mapping $z_0$ to $-z_0$ is just $-g$
due to the form of $f(z)=f_c(z)= z^2+c$. 
\end{notation}

\begin{lemma}[No upper/outside intersections]\label{no-upper}
There is no intersection of $f(\sR_k[2r])$ with
$\Gamma'=f(\Gamma^{\reg}[2r])$ from outside of $\cl f(B(z_0,2r))$.
\end{lemma}

\begin{figure}[h!]

\includegraphics[height=9cm, width=13cm, trim=40pt 50pt 0pt 60pt,clip ]{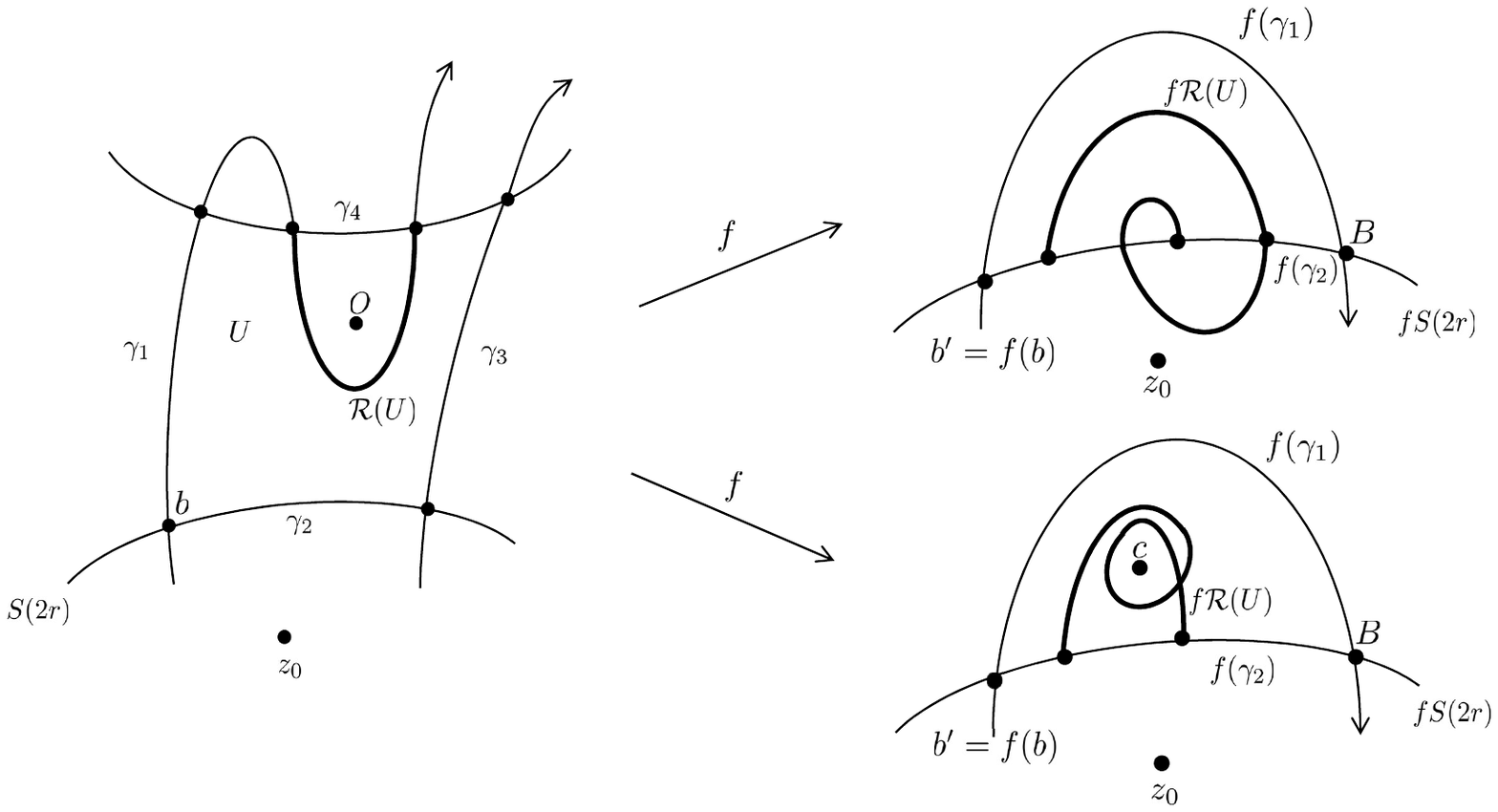}

\caption{By $0\notin  Q^{\reg}[2r]$, a upper intersection would lead to impossible features. }
\label{fig:no-upper}
\end{figure}

\

\begin{proof}
Suppose there is an intersection.
Let $B$ be the first  such an intersection point in $\sR_k[2r]$, except $b'=f(b)$, where $b:=b\sR_k[2r]$.
Then the curves joining $b'$ with $B$ in $\Gamma'$ and in $f(\sR_k[2r])$, denoted $(b',B)$ and
$\sR[b',B]$ respectively, make together a Jordan curve $\delta$, bounding a Jordan domain, containing $c$. The latter follows from the
Fold Lemma \ref{fold}. Then $f^{-1}(\delta)$ forms a Jordan curve, quadrilateral
$\hat\delta$ consisting of four curves

\begin{equation}\label{quadrilateral}
 g(\sR[b',B] ), g((b',B)), -g(\sR[b',B] ), -g((b',B)),
\end{equation}

where $g$ on $\sR[b',B]$ is defined as the extension along $\sR[b',B]$ of $f^{-1}$ starting from $g(b')=b$.
Denote these curves also $\gamma_1, \gamma_2, \gamma_3, \gamma_4$ respectively, the first two in the terminology of the Fold Lemma.
$\hat\delta$ bounds a Jordan domain $U$ containing the point 0.
See Figure \ref{fig:no-upper}

We shall prove now that this is impossible, because $0\notin Q^{\reg}[2r]$ by definition.

The arcs $\gamma_1, \gamma_2$ are in the border $Q^{\reg}[2r]$ and $\gamma_3$ is inside $Q^{\reg}[2r]$ because  close to $\gamma_2$ it is and it cannot intersect $\partial Q^{\reg}[2r]$. The only way 0 occurs outside $Q^{\reg}[2r]$  is to be separated away from it by a curve denoted $\sR(U)$ in $U$ being a part of $\sR_k[2r]$
beginning at entering and ending at exiting $U$, crossing $\gamma_4$ (the upper side of the quadrilateral)
by $\sR(U)$.
Denote the sub-curve of $\gamma_4$ joining the ends of $\sR(U)$ by  $\widehat\gamma_4$

We shall prove that in this situation, $f(\sR_k(U))$ has self-intersections, impossible for external rays. Indeed, the growth of arguments along $\sR_k(U)$ with respect to 0 is less but close to $2\pi$
because $\gamma_4 = -\gamma_2$ is close to $-z_0$, therefore short compared to its distance from 0 in the middle. So the growth of the argument $\Delta$ of $f(\sR_k(U))$ with respect to $c$ is close to $4\pi$.  So $f(\widehat\gamma_4)\cup \sR_k(U))$ must have a self-intersection, otherwise it would be a Jordan curve with $|\Delta|=2\pi$. The only possibility is, that the self-intersection happens in $f(\widehat\gamma_4)$ but then a part of $f(  U )$ is outside of $\delta$ what contradicts the Maximum Principle.

Similarly we prove that for say $k=k_1$ the ray $\sR_{k_2}[2r]$ cannot cross $\gamma_4$ enter $U$
go around 0 and exit through $\gamma_4$.

\end{proof}

\medskip

As used already in \eqref{quadrilateral}, we shall define in the sequel $g$ on $f(\sR_k)$ for $k=k_1,k_2$ as the extension along $\sR_{\s(k)}$ of $f^{-1}$ starting from $g(f(b\sR_k[2r]))=b\sR_k[2r]$. So $g(f(\sR_k))=\sR_k$.
Then for the other branch of $g$ the adequate component of $f^{-1}(f(\sR_k))$ is $-\sR_k$, compare
Notation \ref{not:branches}.
We shall use also the notation $\sR^{\sym}_k$ for $-\sR_k$, or e.g. $\sR^{\sym}_k[2r]$ for restrictions (sym abbreviates \emph{symmetric}).

\

Now let us analyze the intersections of $f(\sR_k[2r]$ with $\Gamma'$ from inside of $f(B(z_0,2r)$,
called \emph{bottom intersection} that we cannot exclude.

\

 Before doing this, consider $\sR_k[r,2r]$. If it intersects $\Gamma$ it must come from $B(z_0,2r)$ going backward from $b\sR_k[2r] = eR_k[r,2r]$ because it does not have another access to $Q^{\reg}[2r]$. Next it leaves $Q^{\reg}[2r]$ also through $\Gamma$ because it must eventually reach $b\sR_k[r]$.

 Denote the points in $\Gamma$ of the entrance and first exit afterward by $K$ and $K'$, the arc in $\Gamma$ joining these two points by
$(K,K')$, the sub-curve joining them in $\sR_k[r,2r]$ by $\Upsilon_K$ in the open domain $Q^{\reg}[2r]$, and Jordan domain between these two curves, by $A_K$. Consider only maximal $A_K$ in the sense of inclusion (compare maximal $A_L$ later on). The exit exists since $\sR_k[r,2r]$ ends (going backward) in $S(r)$.

\medskip

Consider now the curve $-\sR_k[2r]=\sR^{\sym}_k[2r]$. Each its component in $\cl Q^{\reg}[2r]$ must be in one of the sets $A_K$ (or be a point in $\Gamma[2r]$, called a degenerate case, to be discussed later on).  Indeed, the only entrance of this curve to outside $B(z_0,2r))$ is through
the lower arcs, so through $(K,K')$.
It must stay in $Q_0[r]$ if it starts there, because it cannot go around $\sR_k[r]$ and enter e.g. the principal pocket and reach $\Gamma$, intersecting the bridge, because
this is impossible by \eqref{32c}. \underline{See Figure \ref{no-bottom}}.

\begin{figure}[h!]

\includegraphics[height=8cm, width=14cm, trim=20pt 40pt 0pt 20pt,clip]{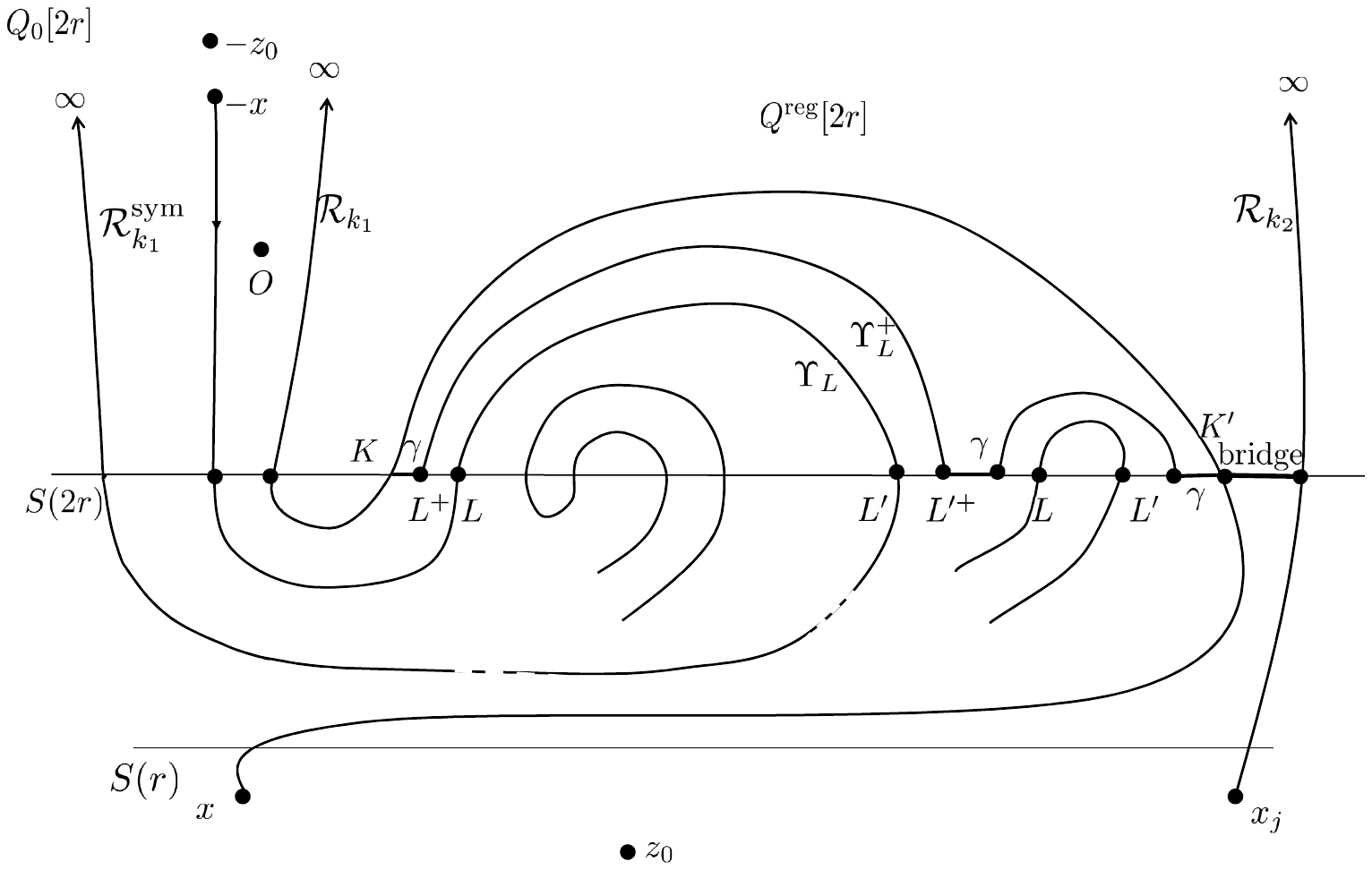}

\caption{Three rays are drawn: $\sR_{k_1}, \sR_{k_2}$ and $\sR_{k_1}^{\sym}$ between them, symmetric to $\sR_{k_1}$, in the singular sector $Q_0$. The critical point is between two latter, in the
quadrilateral as on Figure \ref{fig:no-upper}, but in the critical sector, therefore possible.
The curve $\Upsilon_L$ is maximal and is in the set $A_K$. The corresponding arc $(L+,{L'}^+)$ in the floor $\Gamma$ is pushed beyond $\Upsilon_L$ to $\Upsilon^+_L$ to avoid intersections with $\sR_{k_1}^{\sym}$. Two other curves $\Upsilon_L$ are marked, one of them maximal, so together with its $\Upsilon_L^+$. The ray to which these curves belong is $\sR_{k_1}^{\sym}$, merging at dashed  lines.
}
\label{Gamma}
\end{figure}

\medskip

\begin{figure}

\includegraphics[height=7cm, width=10cm]
{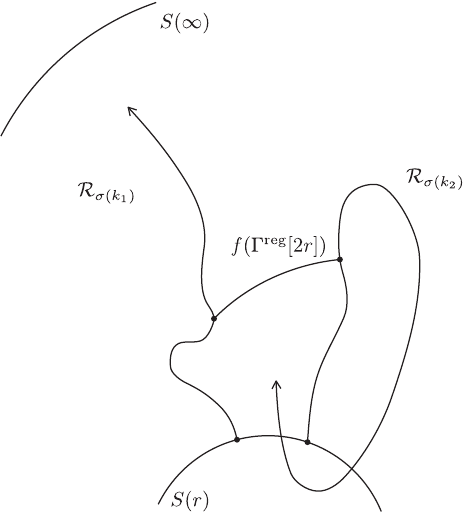}

\caption{No bottom returns. }
\label{no-bottom}
\end{figure}

\

  Given a component of $-\sR_k[2r]$ in $A_K$,
  let us explain the existence of entrance into $A_K$ from $B(z_0,2r)$ through $(K,K')$ at some $L$ and the existence of an exit point  through $(K,K')$  to  $B(z_0,2r)$, at some $L'$, with

  \begin{equation}\label{subarcs}
  (L,L')\subset (K,K').
  \end{equation}

  \underline{The exit exists} because it cannot happen through $\Upsilon_K$ since rays cannot intersect one another and since $-\sR_k[2r]$ escapes to $\infty$, therefore forced to intersect $(K,K')$.

 \underline{The existence of an entrance} follows from the fact that the beginning of  $-\sR_k[2r]$ is $-b\sR_k[2r]\in -S(2r)$
 in $Q_0[2r]$ (the critical sector). Compare Figure \ref{fig:no-upper}, though here
 the quadrilateral containing 0 happens indeed (excluded in $^{\reg}[2r]$. Here $\sR^{\sym}_k[2r]=-\sR_k[2r]$ as above.
 Here $-b\sR_k[2r]$ is in $Q_0[2r]$, so not in $Q^{\reg}[2r]$, into which it must enter therefore through $(K,K')$.


 As above define $\Upsilon_L$ and $A_L$. Since $\Upsilon_L$ cannot intersect one another since they are sub-curves of external rays, the sets $A_L$, hence $\Upsilon_L$, are ordered by the inclusion relation.
Consider only maximal
$\Upsilon_L$ and replace each arc $(L,L')$ for  $k=k_1$ or $k_2$, by a curve $\Upsilon^+_L$ by pushing $\Upsilon_L$ slightly deeper in $A_K$, moving also the ends in $\Gamma$ slightly outside $\cl (L,L')$ to $L^+,L'^+$.

We consider also degenerate cases $\sR^{\sym}_k[2r]$ tangent to $\Gamma$ on the $B(z_0,2r)$ side, where $L=L'$ and $A_L$ and $\Upsilon_L$  consist of single points. ($K$ is irrelevant.) Some can be maximal, some not. For each maximal one we consider $L^+<L=L'<L'^+$.

\medskip

Finally define
\begin{equation}\label{Gamma-star}
\Gamma^*:= \bigcup_L \Upsilon^+_L \cup \gamma,
\end{equation}
where $\gamma=\bigcup \gamma(L)$ is the union of arcs in $\Gamma$ with the ends in
$\{L^+\},\{L'+\}$, making $\Gamma^*$ connected.

\medskip

\begin{remark}\label{both-k}
Note that given $L,L' \in \sR_{k_1}^{\sym}[2r]$ as above, the associated $K,K'$ could be a priori in
any: $\sR_{k_1}[r]$ or in $\sR_{k_2}[r]$ (or given $L,L' \in \sR_{k_2}^{\sym}[2r]$), e.g. in the situation in Figure \ref{channel}.
However it cannot happen that $L_1,L_1'\in \sR_{k_1}^{\sym}[2r]$ and $L_2,L_2'\in \sR_{k_2}^{\sym}[2r]$ for different $L_1$ and $L_2$. This is so because each $L_i, i=1,2$, would involve a different critical value
in $Q_0[2r]$, whereas there is only one: $c$.
A common critical value $c$ would yield an intersection of $\sR_{k_1}[2r]$ with $\sR_{k_2}[2r]$,
because the Jordan domains bounded by arcs of $f(S(2r))\setminus \Gamma'$ and related sub-rays of $\sR_{k_1}[2r]$ and $\sR_{k_2}[2r]$ respectively would both contain $c$.
See Fold Lemma  \ref{fold}, Case 1.
Compare Proof of Lemma \ref{no-upper} and Terminology and Remarks, Item 5, Subsection \ref{Step 2}.
\end{remark}

\

\subsection{Conclusions}\label{Step 6}

By the construction we get

\begin{lemma}\label{f-Gamma}

1. $f(\Gamma^*)$ has no self-intersections.

2. $\Gamma^*$ is disjoint from $-\sR_k[2r]=\sR^{-sym}_k[2r]$, as well as from $\sR_k[2r]$. Hence $f(\Gamma^*)$ does not intersect $f(\sR_k[2r])$ for $k=k_1,k_2$.

3. $f(\Gamma^*)$ is homotopic in $\C\setminus \{c\}$, keeping the ends intact, to
$f(\Gamma^{\reg}[2r])=\Gamma'$.
\end{lemma}

\begin{proof}

1.  $\Gamma^*$ consists of the curves $\Upsilon^+_L$ close to the maximal sub-curves $\Upsilon_L$ of $\sR^{\sym}_k[2r]$ and sub-arcs $\gamma_L$ of $\Gamma^{\reg}[2r]$.


All $f(\Upsilon_L)$ are disjoint from each other, since $\Upsilon_L$'s are
disjoint parts of  $\sR^{\sym}_k[2r]$ on which $f$ is injective. 
So $f$ is also injective on the
the nearby
$\Upsilon_L^+$ 's and  $\Upsilon_L^+$ 's are also pairwise disjoint.
The injectivity also holds on $\gamma$ where $f$ is inverse of $g$.

Notice also that for each $L$, the curve $f(\Upsilon_L)$ does not intersect
$\Gamma'=f(\Gamma)$ except at the ends. Indeed, suppose such a point $B$ of intersection exists.
Then it has at least three $f$-pre-images: one in $\Gamma$, one in $\Upsilon_l\subset \sR^{\sym}_k$,  and one in $\sR_k$. This is impossible since $f$ is quadratic. (Another proof is by using Fold Lemma  \ref{fold}.)

So, referring also to the Maximum Principle,
$f$ is injective on the closure of the Jordan domain bounded by $\Upsilon_L$ and $(L,L')$ for each $L$, in particular the domain does not contain 0. Compare Remark \ref{subfolds}. Hence the injectivity  holds on a small extension of the domain.
Due to being "small" the $f$-images of the
extensions stay pairwise disjoint for different $L$'s .

\smallskip

2. holds because the maximal $\Upsilon_L^+$'s do not  intersect any (not only maximal) $\Upsilon_L$  moreover do not intersect $ \sR^{\sym}_k[2r]$, by construction. Nor $\gamma$
intersects $\sR^{\sym}_k[2r]$.

$\Gamma^*$ is disjoint also from $\sR_k[2r]$ because $\Upsilon_L^+$'s are close to $\Upsilon_L$'s disjoint by definition. Nor $\gamma$ intersects $\sR_k[2r]$ because the former is in $\Gamma \subset S(2r)$ and the latter, except its ends, outside $\cl B(z_0,2r)$ by definition.

  3. holds due to homotopies between the maximal $\Upsilon_L^+$ and $(L^+,L'^+)$
  in the Jordan domain they bound, not containing the point 0. See the item 2. or apply here the Fold Lemma  \ref{fold} for $\gamma_1=(L^+,L'^+)$ and $\gamma_2=\Upsilon_L^+$ or $\gamma_2=-\Upsilon_L^+\subset \sR_k[2r]$.
  These homotopies yield the homotopy we seek for, between the union $\Gamma^*$ in \eqref{Gamma-star} and
  $\Gamma$, by setting the identity on $\gamma$.

  \end{proof}

\begin{corollary}\label{coro}
 For $\partial^*:=\Gamma^* \cup \sR_{k_1}[2r]\cup \sR_{k_2}[2r]$,
 the image: $f(\partial^*)$
 in $\C$ has no self-intersections.
 \end{corollary}

\medskip

Consider $\partial_{\widetilde\C}$ being the boundary of $\widetilde Q^{\reg}[2r]$ in $\widetilde\C$,
 $\partial^+$ being its boundary in $\C$
 and $\partial_\infty$
 being the part of $\partial_{\widetilde\C}$ in $\in S(\infty)$,
 its circle of arguments, where the extension of $f$ is denoted by $F$. Remind that we are in $\widetilde\C$, the plane $\C$ compactified by the circle of arguments $S(\infty)$, where $\widetilde f$ acts as the continuous extension of $f$.
 Namely
 \fbox{$\partial_{\widetilde\C}=\partial^+ \cup \partial_\infty$.}
 Then

 \smallskip

 \begin{lemma}[Main Lemma]\label{injectivity}
 $F$ is injective on $\partial_\infty$, that is the length $|F(\partial_\infty)|< 2\pi$.
 \end{lemma}

 \smallskip

 \begin{proof}

 The growth of arguments
 along the closed curve $\widetilde f (\partial_{\widetilde\C})$ with respect to $c$,
 just the index,  is 0, because $\partial_{\widetilde\C}$ is the boundary of the Jordan domain $\widetilde Q^{\reg}[2r]$,
 which does not contain the critical point 0, nor another $f$-preimage of $c$, because such does not exist, because $f$ is only quadratic and 0 is a double pre-image of $c$.

Define, as in Corollary \ref{coro} and in Subsection \ref{Step 3}, $\partial^*:=\Gamma^* \cup \sR_{k_1}[2r]\cup \sR_{k_2}[2r]$. On $f(\partial^+_\C)$ (the boundary in $\C$), this growth is $\alpha \in (-2\pi, 2\pi)$ because it is the same as on the
homotopic curve $\partial^*$
without self-intersections in the plane. (In fact $\alpha \in (-2\pi, 0)$, because $0\notin Q^{\reg}[2r]$.  So, along  $F(\partial_\infty)$ the growth of argument is $0-\alpha= -\alpha < 2\pi$. The Lemma has been proven.
\end{proof}

\medskip

As a by-product we conclude that the partition of $S(\infty)$ into $\partial_\infty$ and its complement is independent from $r$. Compare Proposition \ref{regular-infty-arc} and
comments in Subsection \ref{Step 2}.

The picture at Figure \ref{impossible}, right hand side, where the length of $F(\partial_\infty)$ exceeds $2\pi$, having a
self-intersection (along an arc) in $S(\infty)$, is impossible.
Otherwise the index with respect to the critical value $\ind_c(\widetilde f (\partial U'))$
 is 1 or 2 (depending on which side of $f(\partial_{\C})$ the critical value $c$ is). But since the only critical point 0
 is not in $U'$, hence $\widetilde\partial U'$ is contractible outside 0, and there is no other pre-image of the critical value $c$ in $U'$ because 0 is its only $f$-pre-image, the index above is 0. A contradiction. 

\medskip

Lemma \ref{injectivity} immediately yields

\medskip

\begin{corollary}\label{coro-order} $F$ preserves the cyclic order of the external arguments of the rays
which begin in $Q_1 \cup ... \cup Q_K$, which is the set of all periodic orbits in $B(z_0,r)$ of minimal period $n$, for $r\le r_0\exp (-n\delta)$, see  Subsection \ref{Step 1}.
\end{corollary}

\medskip

\begin{figure}

\includegraphics[height=6cm, width=12cm]{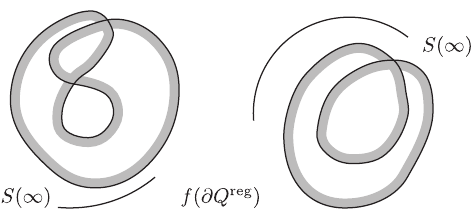}

\caption{Impossible features. On the left hand side the self-intersections in $S(\infty)$ happens, that is $|F (\partial_\infty)| >2\pi$, because of a non-removable self-intersection in $\C$. On the right hand side the marked self-intersection only in $S(\infty)$ is impossible because $0\notin Q^{\reg}[2r]$.}
\label{impossible}
\end{figure}

\

 By the injectivity in Lemma \ref{injectivity},  $F=\widetilde f|_{S(\infty)}$ preserves
 the cyclic order of the boundary (in $S(\infty)$) arcs of the regular  sectors, see Definition \ref{regular-critical}. We can extend it to a homeomorphism $F^o$ so that $F^o$ (combinatorially $\sigma$) maps
 $S(\infty)\setminus \partial_\infty$ to $S(\infty)\setminus F(\partial_\infty)$.
 The latter can be called the \emph{ critical value boundary arc}.

Note, there is no invariant proper sub-cycle of regular sectors under the induced permutation $\sigma$, that is of regular arcs under  $F$.
Otherwise the longest boundary arc in the sub-cycle would be injectively mapped to a twice longer arc, a contradiction. So the permutation $\sigma$ has no cycles of order less than $m$.
So $\sigma$ is just cyclic permutation of order $m$.

Since $\sigma$ maps neighbour sectors (arcs at $\infty$) to neighbour sectors it must be a "rotation" (preserving cyclic order) permutation of our sectors and rays.

In other words $F$ induces an order preserving map on the family of the lifts to $\R$ of the external angles of the rays (sides of the regular sectors, that is all the rays $\sR_j$, and sectors bounded by them,
since $F$ is orientation preserving on $Q^{\reg}(\infty)$ and since $|Q^{\reg}(\infty)|<\pi$, so
"the last (biggest) regular sector" cannot jump over "the first one".

\medskip

The conclusion is that there can be at most 2 orbits of rays in $\bigcup\sR(x,s)$, see \cite[Lemma 2.7]{Milnor}, hence at most $K=2$ periodic orbits of minimal period $n$ in $B(z_0,r)$ 
In fact since the period of $z_0$ is 1 in this Section, there is just one orbit of rays for $n>1$, so $K=1$.

Writing the argument directly, we have proved that $\sigma$ on the sectors, hence on the rays, is transitive, hence has only one orbit. Hence the number $K$ of our periodic orbits of minimal period $n$
is 1. Thus, Theorem \ref{Cremer-fixed} has been proven.
 \end{proof}

\begin{remark} Note that $f$ on each regular $Q_j[r]$ (or $Q_j[2r]$)
is only a local homeomorphism onto its image; it is a priori not clear whether it is injective on $Q_j[r]$ or $Q^{\reg}[r]$. Though the rays cannot intersect one another, the sub-rays $f(\sR_j[r])$ can intersect $f(S(r))$ (if $\sR_j[r]$ pass close to $\{f^{-1}(z_0)\}$.  So we do not know a priori that the order 
of the rays $f(\sR_j[r])$ between their beginnings $f(P_j)$ in $f(S(r))$ and the order of the ends (external arguments) of $\sR_j[r]$ i.e. of
$P_j(\infty)\in S(\infty)$ coincide. Namely it is not clear that such a ray does not meet $S(r)$ again (between $f(P_j)$ and $F(P_j(\infty))$),
i.e. that $f(P_j)$ is its last exit from $f(B(z_0,r))$.

In other words, though the cyclic order of the beginnings $P_j$ of $\sR'_j$ in $S(r)$ coincide with the cyclic order of $f(P_j)$ in $f(S(r))$, it is a priori not clear that $f(P_j)$ are the last exit points from $f(\cl B(z_0,r))$, so it is not clear that their order coincides with the order at $\infty$, as $S(r)$ is not precisely invariant under $f$ (this difficulty is not present if all the rays start at one repelling fixed point, in place of $S(r)$). We resolve it using $Q^{\reg}[2r]$. The rays $\sR'_j$ considered separately could sneak through gaps between
 $f(S(r))$ and $S(r)$.

\end{remark}

\

\subsection{The valence}\label{Step 7}

Here is an easy complementary proposition, used already in Subsection \ref{Step 2}.

\begin{proposition}\label{valence}
For a Cremer fixed point $z_0$ for a quadratic polynomial $f_c$ and $r>0$ small enough
every orbit of period $n > 1$ has the valence $v$ equal at most 2, that is each its point is a limit of at most two external rays.
\end{proposition}
\

\begin{proof}

Let $O=(x_1,...,x_n)$ be a repelling periodic orbit in $B(z_0,r)$, with period at least 2,
Set $\sF=f^{n}$. Let $\sR_1(x_1),...,\sR_v(x_1)$ be all the external rays
starting at $x_1$. Their set will be denoted by $A_1$. Similarly we define $A_j, j=1,...,n$.
Let $Q_s(x_j)$ be sectors between them, that is domains bounded by
$\sR(x_j,s), \sR(x_{j+1 \mod v},s)$ and $\infty$. Then $\sF$ yields $\sigma$  a cyclic (preserving order in external arguments) permutation (combinatorial rotation) of the rays and of the sectors close to $z_0$.

Compare \cite[Lemma 2.1]{Milnor}. See our paragraph [*] in Section 3.

Since $x=x_1$ is fixed for $\sF$ then
$x_2=f(x_1)$ is fixed for $\sF$ as well.
But then, since the sectors $Q_{s,1}$
are pairwise disjoint,
and $\sF$ is a homeomorphism near $B(z_0,r)$, the $\sF$-images of the truncates sectors are contained
in the respective sectors $Q_{\sigma(s),j}$ 
So  $\sigma$ is the identity permutation.

 So for the diagram (phase portrait)
 $A_1,...,A_j$,
 the valence $v$ is 1, or 2 \cite[Lemma 2.7]{Milnor} -- a primitive case.
 See also Appendix.

\end{proof}

See also Lemma \ref{valence2} and Remark \ref{Zakeri}.

\bigskip

\section{Cremer periodic orbits. Proof of Theorem \ref{Cremer-periodic}}\label{Cremer_periodic}

For a Cremer periodic orbit Theorem \ref{Cremer-periodic} holds, generalizing Theorem \ref{Cremer-fixed}.
This section is devoted to this. It will use arguments from Section \ref{Cremer}.

We assume that $r_0$, hence all $r_n$, are small enough that all the closed discs
$\cl B(f^i(z_0),r_n), i=0,...,p-1$ are
pairwise disjoint (far from each other) and far from the critical point 0.
We write $z_i$ for $f^i(z_0)$, for $z_0$ Cremer's of period $p$ and $i=0,...,p-1$.


We can easily choose $z_0$ in the Cremer orbit so that
$|(f^i)'(z_0)|\le 1$ for all $i=0,1,...$. Compare Remark \ref{remark-contraction}.
So the radii of images of $B(z_0,r_n)$ under the differentials $Df^i$ do not grow beyond $r_n$.
So, as in the fixed Cremer point situation, the diameters of $f^i(B(z_0,r_n))$ stay exponentially small for $i\le np$
and distortion stays exponentially close to 1.
We can repeat estimates of Section \ref{distortion}.

\begin{lemma}\label{valence2}
For each orbit of minimal period $np$ for $n>1$, as in Theorem \ref{Cremer-periodic} the valence of all its points (the number of external rays starting at it) is 1 or 2.
\end{lemma}

\begin{proof}
See Proof of Proposition \ref{valence}.
\end{proof}

 Below we shall prove Theorem \ref{Cremer-periodic}. It will follow the proof of Theorem \ref{Cremer-fixed}, with mild complications.
 Originally we attempted  proving Milnor's "pairwise unlinked" condition (4), see Appendix. 
 Unexpectedly proving that condition occurred  difficult (a "proof" in the first version arXiv:2503.03738v1 was incomplete). Here we proceed directly, relying only on the preservation of the cyclic order, Milnor's condition (2), and the proof occurs easy.

 \begin{proof}[Proof of Theorem \ref{Cremer-periodic}]

 Fix an arbitrary $n$. Let $O_1,...,O_K$ be periodic orbits, each of minimal period $np$, entirely contained in the union of the open discs $B(z_i,r)$, for an arbitrary $r\le r_n$. Here
 $r_n=r_0\exp -\delta np$.
 We shall find an upper bound for $K$.



 \

 Consider, as in Section \ref{Cremer}, all external rays $\sR(x,s)$ with their starting points $x\in \bigcup_k O_k $
 with arguments $\theta_s(x), s=0,...m-1$, where $Knp\le m\le 2Knp$. The constant 2 bounds the valence, due to
 Lemma  \ref{valence2}.



 Now, for each ray $\sR(x,s)$ consider $\tau(x,s)$  being the largest $\tau$ so that $\sR(x,s)(\tau)\in \cl B(z_i,r)$, 
 where $x\in B(z_i,r)$.
For each $i$ and $x\in B(z_i,r)$ define the truncated ray $\sR'(x,s))=\sR'(x,s)[r]$, running from $\sR(x,s)(\tau(x,s))$ to $\infty$, as in Section \ref{Cremer}.

\medskip

For each $i$ consider the cyclic order in the family of all the rays $\sR'(x)$ over $x\in B(z_i,r)$ along the circle $\partial B(z_i,r)$ coinciding with the cyclic order of the set $A_i$ of their arguments (angles) in $S(\infty)$,
compare Section \ref{Cremer}. Then the doubling map $F$ restricted to $A_i$ maps it 1-to-1 to $A_{i+1 \mod p}$.
\underline{It preserves the cyclic order}, by the arguments the same as in Proof of Theorem \ref{Cremer-fixed}.
Indeed, we did not use there that $f(z_0)=z_0$.

\medskip

As we already mentioned we would be happy to apply \cite[Lemma 2.7]{Milnor}. 
The sets $A_i$ would play the role of the respective sets
in the formal orbit portrait in \cite[Lemma 2.3]{Milnor}.
For the comfort of the reader we provide the portrait definition in Appendix.
This would give the estimate of the number of periodic orbits of minimal period $np$ being $\exp -\delta n$
close to our Cremer's orbit, at most 2, that is  $K\le 2$.
However we cannot prove
 Milnor's property (4), even for a large $f$-invariant part of the family of rays, namely that all $A_i$'s are
\emph{pairwise unlinked}. Miraculously we do not need it and the proof of Theorem \ref{Cremer-periodic}
using the proof of Theorem \ref{Cremer-fixed} happens to be easy.

\medskip

We use notation from Section \ref{Cremer}. To distinguish objects related to $z_i$ such as e.g.
$\sR_j[r]$ starting in $\partial B(z_i,r)$ 
we write $\sR_j[r]i$, where $i$ is mod $p$. Here $j$ enumerates  our external rays (arguments) starting in $B(z_0,r)$ in the cyclic order. As mentioned above $f$ preserves this order so it is the same cyclic for each $i$. Write $\Gamma_j[r]i$ and $Q_j[r]i$ for the corresponding arcs in
$B(z_i, r)$, defined as in Section \ref{Cremer}, and sectors bounded by $\sR_j[r]i, \sR_{j+1}[r]i$ and  $\Gamma_j[r]i$.

Denote the mapping induced by the doubling map $F$ from $A_i$ to $A_{i+1}$ by $\sigma_i$.
For each $i$ denote the arcs in $S(\infty)$ between cyclic consecutive points of $A_i$ by $S_ji(\infty), j=0,...,m-1$, joining
$\sR_j[r]i$ with $\sR_{j+1 \mod m}[r]i$. Call them \emph{angle arcs}. Each is the $S(\infty)$ part of the boundary of the sector $ Q_j[r]i$. We call the arc (the sector) \emph{regular} if $ Q_j[r]i$ does not contain the critical point 0. For each $i$ the arc (sector)  containing 0 is called \emph{critical}. Assign the index $j=0$ to it. Compare Section \ref{Cremer}.


We shall consider trajectories of the angle arcs, remembering that $F$ is a homeomorphism between
$S_ji(\infty)$ and $S_{\sigma_i(j)}(i+1)(\infty)$ for $S_ji(\infty)$ regular, and makes an additional round along $S(\infty)$ otherwise. For $S_ji(\infty)$ critical there is no choice of the next arc, but to be the arc (gap) $ S(\infty) \setminus\bigcup_{j=1}^{m-1}  F( S_{j}i(\infty)$ by the preservation of the cyclic order. Such a choice of the arc in $F(S_0i(\infty)$ yields the conclusion that the periodic trajectories of the angle arcs are pairwise disjoint.

When iterating $f$, for each $i: 0\le i <p$ the point 0 belongs to just one sector $Q_j[r]i$ for each $i$, that is "serves just one sector". So for the period of our
rays (external arguments)
 $np$ or $2np$ 
 at most $p$ trajectories of the angle arcs can be served.
 Notice that the period is common because (each) $F ^p:A_i\to A_i$ is a combinatorial rotation .

But trajectories consisting of only regular arcs cannot exist because they are periodic so
the length of the arcs doubles at each step tending to $\infty$, in particular exceeding $2\pi$, which is not possible.
So the number of the trajectories of the angles or rays must be upper bounded by $p$.
The trajectories of the (critical) angles cover the whole $S(\infty)$ for each $i$. So,
we obtain the number of their boundary corresponding rays the same, that is $K\le p$.
Compare the end of the proof of Theorem \ref{Cremer-fixed}.


\end{proof}.

\begin{example} The  period 3 \emph{primitive} renormalization example $f(z)=z^2- 7/4$
has orbit portrait at period 3 parabolic orbit
$A_1=\{3/7,4/8\}, A_2=\{6/7,1/7\}, A_3=\{5/7,2/7\}$ with two trajectories of arcs defined above, one meeting a critical arc (sector) once, the other twice, see \cite[Figure 6]{Milnor}. This is the primitive case in Lemma \ref{Lemma-Milnor}.
\end{example}.


\begin{remark} The proof above holds for any polynomial, provided $F:A_i\to A_{i+1}$ preserves cyclic order for each $i$. Then $K\le pd$ where $d$ is the number of critical points in $\C$.
\end{remark}

\begin{remark} Similarly to Theorem \ref{Cremer-periodic} one proves that for a Siegel $S$ disc fixed for quadratic $f$
with the critical point $0\notin \partial S $ and $\partial S$ being a Jordan curve, there are few periodic orbits exponentially close to $S$. The same for any periodic Jordan Siegel disc.
\end{remark}

\section {Hypothesis H. Geometric pressure on periodic orbits}\label{Hypothesis H}

Here we prove

\begin{theorem}\label{H}
Hypothesis H holds for all quadratic polynomials,
with $\sharp\cP\le 2n$.
\end{theorem}

As announced in Introduction, Section \ref{Introduction}, Theorem \ref{H} implies
Theorem \ref{equality-of-pressures} (equality of pressures) due to \cite{PRS2}.

\medskip

The proof of Theorem \ref{H} is similar
 to the one in Section \ref{Cremer} and in a version of a Cremer periodic orbit in
 Section \ref{Cremer_periodic},
 though some additional observations are needed, caused by a possible interference of the critical point. Let us start with two lemmas similar to Lemma \ref{auxiliary} and Corollary \ref{coro1}. They follow in fact \cite[Proposition 3.10]{PRS2} with slightly different proofs; we provide them for completeness.

\medskip
Denote by $\Crit (f)$  the set of all critical points of $f$ in $J(f)$.
For  each $n$ consider the metric
$\rho_n(x,y):=\max_{i=0,...,n-1}\rho (f^i(x),f^i(y))$, where $\rho$ is the spherical metric, see Introduction.

\begin{lemma}[avoiding critical points]\label{Lemma-distortion}
Let $f:\ov\C\to\ov\C$ be a rational function.
Let $\cP$ be a subset of $\Per_n$ as in Hypothesis H with $\rho_n (x,y)< \exp -\d n$, for $n$ large enough and all $x,y\in \cP$. Suppose $\#\cP\ge \exp \d n$.
Then there exists $\epsilon>0$ (not depending on $n$) 
and there exists $\cP'\subset  \cP$ such that

\begin{equation}\label{number}
\#\cP' \ge \exp (n\delta/2),
\end{equation}

and for an arbitrary $x\in \cP'$


\begin{equation}\label{good}
\rho(f^i(x),\Crit(f))\ge \exp (-n\delta/2)
\end{equation}
for all $i$ except at most $a=a(\delta)$ of "bad" ones for a constant $a$, and
for "bad" $i$ and every $ y\in \cP'$

\begin{equation}\label{bad}
\rho (f^i(x), f^i(y))\le (\exp -n\epsilon) \rho(f^i(x),\Crit(f)).
\end{equation}

\end{lemma}

\begin{proof}

There exists $C>0$ depending only on $f$, such that

\noindent $B(\Crit(f), \exp -Cn) \cap \Per_n = \emptyset$. This follows  from
so called Rule I in \cite[Lemma 1]{P-Lyapunov}, saying that for any critical point ${\rm cr}\in J(f)$ and $x, f^N(x)\in B({\rm cr},r)$ it holds $N > C^{-1}\log 1/r$, applied to $n=N$ and $r=\exp -n$.

\medskip


Denote the critical points in $J(f)$ by $0_1,...,0_D$.
Consider the annulus $A_j=B(0_j,\exp (-n\delta/2))\setminus B(0,\exp (-Cn))$.
In the logarithmic coordinates $A_j$ (more precisely $\log (A_j-0_j)$, with the branch of logarithm with range in $\{x+iy: -\pi<y\le \pi\}$),
can be covered, taking in account all $0_j$, by
$$
D \frac{2\pi (Cn - n\delta/2)}{\exp -2n\epsilon} \le \exp 3\epsilon n
$$
squares (in the logarithmic coordinates) of the sides of length $\exp -n\epsilon$ for $n$ large enough.
Enumerate them by $Q_s$.

If $i$ is bad for $x\in\cP$, then $f^i(x)$ belongs to one of these squares for some $j$.
We select $s=s(i)$ with the largest number of points in $f^i(\cP)$ belonging to it, so
$$
s(i)\ge \frac{\exp \delta n}{ \exp 3\epsilon n}  = \exp n (\delta - 3\epsilon).
$$
This (first) selection
defines
$$
\cP'_1: =\{y\in \cP: f^{s(i)}(y)\in Q_s\}.
$$

 Next, among other bad $i$'s we make a selection in $\cP'_1$ and obtain  $\cP'_2$, etc. until $\cP'_a$ having at least
$$
\exp(n(\delta-3Da\epsilon))
$$
elements for $n$ large enough.

Here $a$ bounds from above the number of returns of the trajectory

\noindent $f^i(x)_{i=1,...,n-1}$ to
$B(0_j,\exp(n\delta/2))$ for each $j$. It is inverse of the smallest gap between consecutive times of bad returns near $0_j$, which by Rule I is bounded below by
$C^{-1}n\delta/2$. Thus
$a\le 2C/\delta$.

We conclude \eqref{number}, by taking any $\epsilon <  \delta^2/12DC$. The inequalities \eqref{good} and \eqref{bad} follow immediately from the definitions.

\end{proof}

 \begin{lemma}[small distortion]\label{auxiliary2}
 Let $f:\ov\C\to\ov\C$ be a rational function. Then for every $\delta>0$ there exist $\epsilon>0$
 and $\beta>0$               
 such that for all        
 $n$ large enough the following holds. Let
 $f^i(x), f^i(y), i=0,..., n-1$ be two different $f$-periodic orbits of minimal period $n$ which are $r_n = \exp -n\delta$ close in the metric $\rho_n$, that is such that
 $\rho_n (x,y) \le r_n$. Set $\rho(x,y):=r\le r_n$.

Assume  that the conditions either \eqref{good} or \eqref{bad}
are satisfied for each $i$, as in Lemma \ref{Lemma-distortion}. 
Then

\medskip

1. All $f^i(B(x,(1+\beta)^n r)), i=0,...,n-1$ are uniformly boundedly distorted, with distortion of order at most $1+\exp(-n\epsilon)$ and

 \begin{equation}\label{modul}
2. \  \  \hbox{the distance} \ \ \dist(f^i(B(x,(1+\beta)^n r)),\Crit(f)) \ge
\end{equation}
$$
\exp (n\epsilon) \diam f^i(B(x,(1+\beta)^n r)).
$$

 In particular $f^i$ is injective on $B(x,(1+\beta)^n r)$ for each $i=0,...,n$.

 \end{lemma}

Geometrically, $f^i$ is almost a similarity on $B(x,(1+\beta)^n r)$ for each $i=0,...,n$.

 \begin{proof}  

Write $x=x_0, y=y_0.$  
Notice that due to \eqref{good} and \eqref{bad} we can neglect the influence  of critical points. So  
\begin{equation}\label{ratio}
\Big| \frac
{|\rho(f^{i+1}(x_0), f^{i+1}(y_0))/\rho(f^i(x_0), f^i(y_0))}{|f'(f^i(x_0))|} - 1\Big| \le
\exp -n\epsilon.
\end{equation}
compare \eqref{complex distortion}.

 To get this, repeat the proofs of Lemma \ref{auxiliary} and Corollary \ref{coro1} to estimate $f'(x_i)$ consecutively and composing them. The same for every other $y\in B(x,(1+\beta)^n r)$. This yields bounds on distortion.

 The orbit $(x=x_0,...,x_{n-1},x_n=x_0)$ plays the role of Cremer's periodic $(z_0,...,z_n=z_0)$. Here we consider $f^i(B(x,(1+\beta)^n r)$ in place of $f^i(B(z_0,r))$.
 Note that here $|f'(f^i(x))|$ need not be 1, so $|(f^i)'(x)|$ can grow exponentially for a while,
 though $(f^n)'(x)\approx 1$ by \eqref{ratio} and the periodicity of $x$ and $y$ of period $n$.
 The sets $f^i(B(x,(1+\beta)^n  r))$ stay exponentially small since the diameters $\diam f^i(B,r)$ do not exceed
 much $\rho_n(x,y) \le \exp -\delta n $. Thus the nonlinearity of $f$ along $f^i(B(x,(1+\beta)^n r))$ stays exponentially small.

\end{proof}

\

\begin{remark}\label{remark-contraction}
 In the periodic orbit of $x$ we can choose a point, say $x_0$, such that
\begin{equation}\label{contraction}
|(f^i)'(x_0)|\le 1 \ \ \hbox{for all}\ \ i=0,1,..., n-1.
 \end{equation}
 Indeed, choose $i$ for which
 $|(f^i)'(x)|$
 attains maximum and define $x_0:=f^i(x)$.

 \medskip

 Nevertheless we did not need  to apply this choice in the proof of Lemma \ref{auxiliary2}, because $\diam f^i(B(x,(1+\beta)^n r))$ was
 controlled to be small because of their distance from critical points. The added margin $(1+\beta)^n$
 stays controlled under iteration as in Section \ref{distortion}.

 We shall apply Remark \ref{remark-contraction} in the proof of Theorem \ref{H}, below.
 \end{remark}

 \bigskip


In the proof of Theorem \ref{H} we shall refer  to Theorems \ref{Cremer-fixed} and Theorem \ref{Cremer-periodic}, but not directly,
because those theorems treat the situations in presence of Cremer fixed point (or periodic orbit).

Instead we shall apply the following fact having the same proof:

\begin{theorem}\label{non-periodic} For every quadratic polynomial $f=f_c$, for every $\delta>0$ there exist $C>0$ and $r_0>0$ such that for every integers $p$ and $n$ and every sequence of discs $B_i, 
i=0,...,p-1$ centered in $J(f)$, if

 $\bullet$\ \ $\diam B_i \le r_0\exp(-\delta np) \times \dist(C B_i,\Crit(f))$,

 \noindent where for $B=B(x,r)\subset \C$ we write $CB:=\{y: \rho(x,y) < Cr\}$,

 $\bullet$\ \  $\diam B_i \le r_0\exp(-\delta np)$,

 $\bullet$\ \ and $C B_i$ are pairwise disjoint,

then there are at most $p$ points $x\in B_0$ of minimal period $np$,
such that $f^{kp+i}(x)\in B_i$ for each $k:0\le k < n$ and $i:0\le i < p$.
\end{theorem}

The proof is the same as in presence of a Cremer orbit of minimal period $p$, whose role is taken over by
any periodic orbit of $x$ in its statement, of minimal period $np$. It does not matter whether as in Section \ref{Cremer}
the sets $B(f^i(z_0),r)), i=0,1,...,p-1$ are considered, or $B_i$ above. We can assume $B_i$ are centered at $f^i(x)$.
The maps $f^i$ are almost conformal affine on $B_0$ due the bulleted assumptions and Section \ref{distortion}.
Preservation by $f$ of order of the external rays involved, between $B_i$ and $B_{i+1}$, can be proved similarly as in Sections \ref{Cremer} and \ref{Cremer_periodic}.

\

\begin{figure}

\includegraphics[height=6cm, width=12cm]{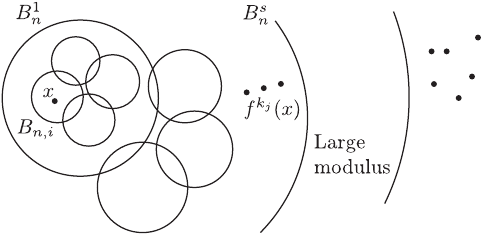}

\caption{Moduli}
\label{moduli}
\end{figure}

\begin{proof}[Proof of Theorem \ref{H}]

Fixed $n$ given $x\in \cP'$ 
choose $x_0$ in the periodic orbit of $x$
such that $|(f^i)'(x_0)|\le 1$ for each $i$, according to Remark \ref{remark-contraction}.
Consider the disc $B_n:=B(x, r)$ for
$r=2\max\{\rho(x,y): y\in\cP'\} \le r_n=r_0 \exp-\delta n$
as in
 Lemma \ref{auxiliary2} (we can omit $r_0$ for $n$ large enough). So all $f^i, i=0,...,n$ are almost similarities on $B(x,(1+\beta)^n r)$, compare notation in Theorem \ref{Cremer-periodic}. Here $\cP'$ is as in Lemma \ref{Lemma-distortion}, in particular satisfying \eqref{number}.

\noindent Write
$B_{n,i}:=f^i(B_n)$.

We take $C>1$ large enough that the proof below makes sense.
Let $i=k_1$ be the first time of return  of $x$ under the action by $f^i$ to $CB_n=B(x,Cr)$.
Since $f^n(x)=x\in CB_n$, $k_1\le n$.
There are two cases now

\

I. $k_1=n$.
Then we have the situation similar to that in Theorem \ref{Cremer-periodic}.
The quasi-discs
$f^i(B(x, \frac14 Cr))$ for all $0<i<n$
are disjoint from $B(x,\frac14 Cr)$.
In particular $B_{n,i}$ are well disjoint from $B_n$ for $C$ large enough.

Hence, for $0\le i < n$, \
$B_{n,i}$
are pairwise disjoint, with annuli separating each of them from the others
of modulus at least $\frac1{2\pi}\log C/3$. We applied here Remark \ref{remark-contraction}.
The pairwise disjointness follows, since otherwise, for some $0<i<k<n$, we would  have
$f^{n-i}(B_{n,i})\cap f^{n-i}(B_{n,k})\not=\emptyset$.
The latter quasi-discs are contained in $B(x, 2r)$ and $f^{k-i}(B(x, r))$
which therefore intersect. Contradiction. The same consideration holds for moduli.

\medskip

Finally we apply Theorem \ref{Cremer-periodic} with $p=n$.
So the number of periodic orbits of minimal period $n$ in $\cC'$ we have discussed, does not exceed $n$. This contradicts \eqref{number}.

Note that no Cremer periodic orbit is involved; it is not needed. See Theorem \ref{non-periodic}.

\

II. $k_1< n$.

\

See Figure \ref{moduli}.


\

Consider $\bigcup_{i=0}^{n/k_1} f^{ik_1}(B(x, Cr))$ and the disc $B_n^1:=B(x, n4Cr)$ containing it.
This containing follows from Remark \ref{remark-contraction}, for a right choice of $x$ in its orbit.
Let $i=k_2\le k_1$ be the first time $i$ of return  of $x$ under the action by $f^i$ to $CB_n^1$.
Then either

$I_2.$\, $k_2=k_1$. Then analogously to the case I. we consider the pairwise disjoint quasi-discs
$B_{n,i}^1:= f^i(B_n^1)$ for $i=0,...,k_1-1$ and we consider $p=k_1$

or

$II_2.$\, $k_2<k_1$. Then we define $B_n^2$ and $B_{n,i}^2$ for $i=0,...,k_2-1$.

Etc. until $k_s=k_{s-1}$ or $k_s=1$ for certain $s$.

Notice that $s\le \log_2 n$, since for each time $i<s$ we have $k_{i+1} \le k_i/2$.

Notice that due to our construction $\diam B_{n,i}^s \le r (4C)^s \prod_{i=0}^{s-1} k_i/k_{i+1}$, where $k_0:=n$.
So, for each $i: 0\le i < n$,
\begin{equation}\label{s-step}
\diam B_{n,i}^s \le r (4C)^{\log_2 n} n  \le r n^{\log_2 4C },
\end{equation}
where $n^{\log_2 4C }$ grows slower than exponentially with $n$, in particular slower  than $(1+\beta)^n$.

So, if $k_s=1$ then we have the situation as in Section \ref{Cremer} and $B(z_0,r)$ and $r$ exponentially small.
If $k_s=k_{s-1}$ we have the situation as in Section \ref{Cremer_periodic} with $p=k_s$.
Again we receive the number of periodic orbits we have discussed does not exceeding $p\le n$.
This contradicts \eqref{number}.

\end{proof}

\

\section{Remarks by C.L. Petersen and by S. Zakeri}

\

\begin{remark}\label{baby-M}
Theorem \ref{Cremer-periodic} can be strengthened, replacing at most $p$ orbits of minimal period $np$,
by at most one orbit, for all $n>1$. This has been told to me by Carsten Lunde Petersen.
Indeed, one can conclude this from Theorem \ref{Cremer-fixed} via renormalizations. See \cite{DH} and
\cite[Theorem 1 de Modulation]{Haissinsky}.

First notice that the Cremer orbit parameter $c$, and more generally such that $f=f_c$ has an indifferent cycle  $O_f(z_0)$ with multiplier $\lambda\not=1$,  belongs to the boundary of a hyperbolic component $M'$ of the Mandelbrot set $M$, where the corresponding cycle has period $p$ and is indifferent. This is a classical fact proven with the use of Implicit Function Theorem. It is sufficient to consider only a primitive case, where the multiplier is different from 1. Then the existence of a quadratic-like couple $f_c^p:U'_c\to U_c$ follows from \cite{DH, Haissinsky}. $f_c^p|_{U'_c}$ is quasi-conformally conjugate to a restriction of a quadratic polynomial $f_{\chi(c)}$. So by topological arguments the conjugacy preserves the multiplier at $z_0$ and by H\"older property preserves scales up to $\delta$ in
$\exp -\delta n$. In conclusion Theorem \ref{Cremer-fixed} yields the strengthened Theorem \ref{Cremer-periodic}.

Indeed, for any $f$-orbit of period $np$ close to $O(z_0)$, its $f^p$-suborbit is close to $z_0$. So as not escaping from the small Julia set $J'\ni z_0$ for $f_c^p|_{U'_c}$, it must belong to $J'$ and by
Theorem \ref{Cremer-fixed} there is at most one such orbit.

\end{remark}
\begin{remark}\label{Zakeri}

Saeed Zakeri has informed me that Proposition \ref{valence} can have versions much stronger,
following from his \cite[Theorem 5]{Zakeri}. E.g. for quadratic polynomials with Cremer fixed point $z_0$, there is no need to assume that the orbit (cycle) of period $n$ is small and that in such general circumstances there are no bi-accessible points in Julia set at all, except maybe those whose forward trajectories hit $z_0$.

\end{remark}

\

\section{Appendix. Orbit portraits, following \cite{Milnor}}\label{App-Milnor}

\

Let $A_1,...,A_p$ be subsets of the circle of arguments (angles) at $\infty$.
This collection is called \emph{the formal orbit portrait} if the following conditions are satisfied:

\medskip

(1) Each $A_i$ is a finite subset of $\R/2\pi\Z$.

(2) For each $i \mod p$, the doubling map $F: \theta \to  2\theta (\mod 2\pi\Z )$ carries $A_i$ bijectively
onto $A_{i+1}$ preserving cyclic order around the circle,

(3) All of the angles in $A_1\cup ... \cup A_p$ are periodic under doubling, with a common
period $r p$, and

(4) the sets $A_1, ..., A_p$ are pairwise unlinked; that is, for each $i\not=k$ the sets $A_i$
and $A_k$ are contained in disjoint sub-intervals of $\R / 2\pi\Z$.

\


\begin{figure}[h]

\includegraphics[height=6cm, width=9cm]{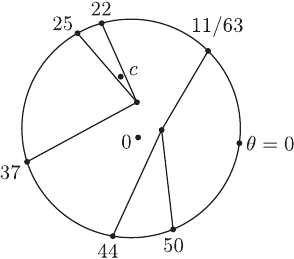}

\caption{Formal orbit portrait, \cite[Fig. 2]{Milnor}}
\label{portrait}
\end{figure}

\

Here at Figure \ref{portrait}, is Milnor's schematic diagram illustrating the formal orbit portrait associated to $z\mapsto z^2+c$ for $c=\frac14 e^{2\pi i/3}-1$ for a period 2 parabolic orbit, with its external ray orbit of period 6.

\noindent $A_1/2\pi = \{22/63, 25/63, 37/63\}$ and $A_2/2\pi = \{11/63, 44/63, 50/63\}$. The critical value $c$ lies in the smallest sector. Compare \cite[Figures 1 and 2]{Milnor}.

\

\begin{lemma}[Milnor, Lemma 2.7, primitive vs satellite]\label{Lemma-Milnor}

Any formal orbit portrait of valence $v>r$ must have $v=2$ and $r=1$.  It follows then that there are just two possibilities:

Primitive Case. If $r=1$ so that every every ray which lands on the period $p$ orbit is mapped by to itself by
$f^p$, then at most two rays land on each orbit point.

Satellite Case. If $r>1$, then $v=r$, so that exactly $r$ rays land on each orbit point, and all of
these rays belong to a single cyclic orbit.

\end{lemma}

The number $r$ is defined as the number of points in each $A_i$ in one $F$-orbit (by (3) it does not depend on the orbit nor $i$).


\

\

\end{document}